\documentclass[sn-mathphys,Numbered]{sn-jnl}

\pdfoutput=1

\usepackage{graphicx}%
\usepackage{multirow}%
\usepackage{amsmath,amssymb,amsfonts}%
\usepackage{amsthm}%
\usepackage{mathrsfs}%
\usepackage[title]{appendix}%
\usepackage{xcolor}%
\usepackage{textcomp}%
\usepackage{manyfoot}%
\usepackage{booktabs}%
\usepackage{algorithm}%
\usepackage{algorithmicx}%
\usepackage{algpseudocode}%
\usepackage{listings}%

\usepackage[utf8]{inputenc}
\usepackage[english]{babel}
\usepackage[T1]{fontenc}
\usepackage{amsmath}
\usepackage{hyperref}

\usepackage{doi}

\usepackage{graphicx}

\usepackage[normalem]{ulem}

\usepackage{amssymb}

\usepackage{dcpic,pictexwd}

\usepackage{amsfonts}
\usepackage{graphicx}
\usepackage{epstopdf}

\newcommand{\intprod}{\mathbin{\raisebox{\depth}{\scalebox{1}[2]{$\lrcorner$}}}}

\newtheorem{mydef}{Definition}
\newtheorem{mytheorem}{Theorem}
\newtheorem{mylemma}{Lemma}
\newtheorem*{myremark}{Remark}
\newtheorem{myrp}{RP}

\newcommand{\rmd}{{\rm d}}
\newcommand{\TM}{{\mathrm{T}\Omega}}
\newcommand{\TxM}{{\mathrm{T}_x\Omega}}

\newcommand{\TstarxM}{{\mathrm{T}^{*}_x\Omega}}

\begin{document}

\title[Article Title]{An electrical engineering perspective on naturality in computational physics}


\author*[1]{\fnm{P. Robert} \sur{Kotiuga}}\email{prk@bu.edu}

\author*[2]{\fnm{Valtteri} \sur{Lahtinen}}\email{valtteri.lahtinen@quanscient.com}


\affil[1]{\orgdiv{Department of Electrical and Computer Engineering}, \orgname{Boston University}, \orgaddress{\street{8 St. Mary's St}, \city{Boston}, \postcode{02215}, \state{MA}, \country{USA}}}

\affil[2]{\orgname{Quanscient Oy}, \orgaddress{\street{Åkerlundinkatu 8}, \city{Tampere}, \postcode{33100}, \country{Finland}}}

\abstract{We look at computational physics from an electrical engineering perspective and suggest that several concepts of mathematics, not so well-established in computational physics literature, present themselves as opportunities in the field. We discuss elliptic complexes and highlight the category theoretical background and its role as a unifying language between algebraic topology, differential geometry, and modelling software design. In particular, the ubiquitous concept of naturality is central. Natural differential operators have functorial analogues on the cochains of triangulated manifolds. In order to establish this correspondence, we derive formulas involving simplices and barycentric coordinates, defining discrete vector fields and a discrete Lie derivative as a result of a discrete analogue of Cartan's magic formula. This theorem is the main mathematical result of the paper.}

\keywords{Finite element modelling, Discretization, Computational physics, Applied mathematics}



\maketitle

{
	
\begin{center} \it For Alain Bossavit, with respect for his scientific humility and mentorship. \end{center}

}

\section{Introduction}

A \emph{finite formulation} of the governing equations is necessary for computer simulations of physical phenomena. The traditional approach to this in continuum physics is a discretization -- a way to transfer from continuum to discrete -- for the spaces arising from the (partial) differential equations that describe the phenomenon in question. In this paper, we discuss finite formulations from an electrical engineering perspective, reaching out especially to the engineering, computational physics, and applied mathematics communities. The topic of our interest is the nature of the process of obtaining finite, discrete descriptions of (continuum) physics. Through building on the foundations of geometric, global analysis \citep{Eells}, \citep{Palais1968} and articulating physics via such concepts as natural differential operators and elliptic complexes, the framework of \emph{category theory} \citep{Adamek}, \citep{Riehl2017} provides a sense of synthesis to the seemingly scattered field. At first glance, the concepts of category theory may steer the reader's thoughts towards foundations of mathematics. However, on a closer look, it is also the glue and common ground throughout computational physics. Evolving from such foundational questions as the distinction between sets and classes in order to resolve the paradoxes of set theory and the foundations of mathematics, category theory leads to the foundations of computer science. Unfortunately, these days, the foundations of computer science are largely eclipsed by the needs of computer engineering. However, computational electromagnetics has demonstrated a more direct role for category theory in linking engineering design to software methodology, as we will discuss in this paper.

Inspired by Freeman Dyson's paper entitled \emph{Missed opportunities} \citep{Dyson}, throughout this paper, we aim to bring forth concepts of mathematics, which, at least in hindsight, present themselves as opportunities for computational physics but paths less traveled in this context. Certainly, we are not in the position to give such insightful advice as Dyson did in his aforementioned influential work, but we hope to be able to frame some future research with a fresh point of view stemming from our electrical engineering background. Hence, this paper is not about historically missed opportunities, but a \emph{programmatic} paper, which through reflecting history and framing future, suggests a research program spanning pure and applied mathematics, computational physics and engineering.

\subsection{Background and motivation}

From an engineer's viewpoint, the importance and practical value of considering foundational issues related to discrete descriptions of continuum physics is evident in the context of framing forward-looking research problems. In electrical engineering, albeit via arguably confusing terminology, Gabriel Kron was one of the pioneers to promote the use of both differential geometry and algebraic topology in applications, and for finding unifying concepts in discrete and continuous descriptions of electromagnetism (see e.g. \citep{Kron1959a}, \citep{Kron1959b}). J. P. Roth would also follow up on Kron's work, and discuss e.g. homological techniques for electrical engineering in a systematic manner \citep{Roth1955}, \citep{Roth1959}, \citep{Roth1971}. Later, \citeauthor{Balasubramanian} explained Kron's methods with more modern terms partially filling the gap between the languages used by Kron and mathematicians \citep{Balasubramanian}.

Taking a view from the 1980s in the computational electromagnetics (CEM) community, the ideas initiated by Kron and his followers had still not quite penetrated the field. Finite element modelling was becoming more accessible due to development of computers. However, for example issues regarding spurious, unphysical modes had plagued modelling in the finite elements setting since the 1960s, when Silvester brought finite elements to the context of CEM \citep{Silvester}. Meanwhile, methods of finite-difference time-domain (FDTD) type utilizing staggered grids, originally introduced by Yee (\citeyear{Yee}), were also being developed \citep{Taflove}. In the wake of recent developments by mathematicians such as Dodziuk (\citeyear{Dodziuk1976}), (\citeyear{Dodziuk1981}), Kotiuga's works of the time (see e.g. \citep{Kotiuga1984}, \citep{Kotiuga1989}) presented a journey to the early mathematical literature, while connecting with and renewing the state-of-the-art of CEM, emphasizing topological aspects. Around the same time, Alain Bossavit described the framework of \emph{Whitney forms} for finite element setting, introduced by Whitney long before the era of computer-aided modelling \citep{Whitney}, to the engineering community \citep{Bossavit}, \citep{Bossavit1988}, \citep{BossavitWhitneyComplex}. After the appearance of Kotiuga's PhD thesis \citep{Kotiuga1984}, Bossavit noted the correspondence between Whitney forms and N\'ed\'elec elements \citep{Nedelec} and their affine invariance, and that this was a key to resolving several outstanding problems in three-dimensional computational electromagnetics. In fact, Whitney forms first appeared in \citep{Weil}. This paper by Andr\'e Weil, representing partially a foundation of algebraic topology and a category theoretical rationale for Whitney forms, was a source of considerable anxiety, initiating the final countdown for Bourbaki's aspirations for a set-theoretic foundation of mathematics. In the meantime Whitney forms played a crucial role in some geometrically motivated mathematical works by, for example, Dodziuk. However, it was the engineering community that finally revived Whitney forms in the 1980s, understanding their significance in the context of CEM. For engineering applications, Whitney forms, widely used today in the context of finite element computations, ensure structure-preservation, continuity properties, and help avoid unphysical, spurious solutions \citep{Bossavit1990}, thus redeeming some of the issues that were related to finite element methods. Later, Ralf Hiptmair re-formulated and generalized the framework to higher-order forms \citep{Hiptmair1999}, \citep{Hiptmair2001a}.

While the framework of Whitney forms discretizes metric-independent properties effectively, on the flip-side of the coin one has the problem of discrete constitutive relations, which tend to be metric-dependent \citep{Tarhasaari}, \citep{Hiptmair2001b}, \citep{Kangas}, \citep{Auchmann}. As the search for compatible, structure-preserving discretizations led to the research programs of \emph{discrete exterior calculus} (DEC) \citep{Hirani} and \emph{finite element exterior calculus} (FEEC) \citep{Arnold2006, Arnold}, Kotiuga brought some of the relevant historical mathematical questions in relation to such program together with present state of CEM \citep{Kotiuga2008}. In particular, the \emph{discrete star localization problem} (DSLP), which states that a discrete version of Hodge star operator, ubiquitous in constitutive laws of physics, cannot be local in the same sense as its continuum counterpart it mimics, casts a shadow on the hope of discretizing constitutive relations in a fully structure-preserving manner \citep{Bossavit1999, Tarhasaari}. Moreover, the so-called \emph{commutative cochain problem}, related to the impossibility of finding a fully structure-preserving discrete analogue for the exterior product of differential forms, poses further limitations for the possibilities of DEC and FEEC. Simply, some of the properties of continuum theory are not reproduced in the discrete setting, and this is reflected in the modelling decisions one has to make. The bottom line is that the inherent difference of discrete and continuum formulations is thus still largely unresolved. This sets one out to look for a coherent theoretical framework (relevant categories and functors) rendering a given discretization procedure as canonical (functorial) within it.

As for continuum descriptions, theories of physics typically formalize our intuition of continuum by stating their principles in the form of differential equations required to hold in each point of the space. As a consequence, one is confronted with infinite-dimensional function spaces when utilizing such theories to make predictions. Apart from very simple cases where closed-form solutions are possible, approximate solutions to the arising field problems are needed: One is required to form a discretization for the function space to find a solution with a finite set of information. Numerical solution methods of field problems, such as the finite element method (FEM)\footnote{Differential geometry provides a modern framework and a suitable formalism for presenting FEM. The historical roots go back to ``Dirichlet's Principle'', which is one of Hilbert's 23 problems from 1900. Historically, it is not due to Dirichlet but Riemann who attributed the principle to Dirichlet because he needed to assume the existence of solutions to Laplace's equation in order to develop his theory of functions of a complex variable. Hilbert's student, Ritz, attacked this problem from a constructive point of view in the early 1900s, but it was Richard Courant's definitive treatment in the late 1920s that gave birth to both a constructive proof of Dirichlet's principle, and to FEM. As a numerical method, FEM can be traced to the interactions between John Von Neumann and Courant close to two decades later in the early days of programmable computers. For the related history, see e.g. \citep{Courant}, \citep{Gander}, \citep{Gorkin} \citep{Pelosi}, \citep{Taylor}, and references therein.} \citep{Bossavit1998b} \citep{Oten}, \citep{Reddy}, \citep{Strang}, rely on such discretizations; they provide a mapping from continuum to a discrete space.

Modelling the world directly utilizing macroscopic quantities, in some sense, avoids the tedious limit process that is inevitable in the differential description of the world \citep{Tonti2014}. Moreover, one can argue that macroscopic quantities are the ones closer to our true experience, since they are the ones closer to measurements: Instead of electric fields we measure voltages, and instead of electric current densities we measure electric currents, for example.\footnote{Of course, this is not to say that we are claiming anything about \emph{direct} measurability of any macroscopic quantity.} Numerical methods starting from algebraic equations concerning macroscopic quantities, such as the cell method \citep{Tonti2001}, represent this point of view.\footnote{The fundamental ideas leading to the cell method date back to Maxwell's notions of analogies between different physical theories and were articulated by Tonti in terms of Tonti diagrams \citep{Tonti1972}, \citep{Tonti1977}. Continuing this path, now applied category theory can be seen as a pervasive effort of transporting analogous ideas formally between mathematics and science as well as between different fields of science. For some striking analogies articulated via categories, see e.g. \citep{Baez}.} The cell method does not present a model in terms of discretized differential equations to obtain a finite formulation suitable for computers, rather requiring \emph{algebraic} equations to hold for all geometric entities of relevance in a continuous space before discretization. Methods such as the cell method, or the finite integration technique (FIT) \citep{Weiland}\footnote{Even though FIT explicitly considers physical variables as integrals of \emph{field} quantities, one is still effectively dealing with the physics macroscopically to begin with. Hence, we consider FIT a close relative of the cell method. FIT indeed started off dubbed as a FDTD-like method but was later explicitly set in a more geometric setting by Weiland \citep{Weiland96}.} are sometimes called \emph{directly discrete methods} as they would \emph{seem} to provide a discrete formulation in a more direct sense than FEM. Ultimately, however, such discretization methods lack a rigorous framework for error analysis.

\subsection{Methodology and objectives}

For the sake of simplicity, we will be concerned with representing physics with partial differential equations expressed in terms of \emph{differential forms} and algebraic equations expressed in terms of \emph{cochains}, a macroscopic counterpart for differential forms. Thus, mixed tensor physics is not explicitly considered. 
Cochains are the mathematical tool for representing macroscopic quantities of physics. There is a well-known  correspondence between cochains and differential forms \cite{Whitney}. Through limit process, one obtains their differential form counterparts, defined pointwise, from which through integration over macroscopic $p$-volumes, one again obtains cochains. For example, in electromagnetism, the integral of the magnetic flux density differential 2-form $B$ over a surface $S$ yields the magnetic flux cochain $\Phi$ through $S$. Knowing the magnetic fluxes through \emph{all} surfaces is the same as knowing $B$. We emphasize the similarity of different routes to a discrete formulation.\footnote{For comprehensive studies of similarities between different numerical methods from a different viewpoint, see e.g. \citep{Bochev}, \citep{Mattiussi}.} In particular, we discuss Galerkin FEM and cochain methods as examples. While cochain methods access the global macroscopic quantities directly unlike FEM, infinite-dimensional function spaces are inherently present in both before the discretization is made. This has implications for computational physics through the identification of structural similarities.

We will highlight the category theoretical background arising from algebraic topology and differential geometry within computational physics. In particular, the ubiquitous concept of \emph{naturality} plays a central role \citep{Kolar}. There are many facets to such background. Category theory has been slowly penetrating modern physics and computation \citep{Categorification}, \citep{Baez}, \citep{Coecke}, \citep{Lal}, \citep{Nikolaus} and even engineering sciences \citep{CategoriesInControl}, \citep{PassiveLinear}. Semantics of computation may be articulated via category theory: Lambda calculus and functional programming sit naturally in this context \citep{Baez}, \citep{Milewski}. Common functional programming languages, such as Haskell\footnote{\url{http://www.haskell.org}}, can be seen as instances of the general theory. As category theory unifies concepts in algebraic topology, differential geometry, and computer science, we aim to argue that category theory may provide synthesis to the seemingly scattered field of computational physics and be beneficial in programming discrete solvers.

Hence, this is a paper about computational physics, written from the viewpoint of electrical engineers and meant to initiate discussion across disciplines. As recent engineering history shows, rigorous identification of underlying structures is often the key for forward-looking research. As a result, we frame and suggest a research program.

\subsection{The structure of the article}

In section~2, we discuss mathematical structures underlying continuum physics and its discrete descriptions. In particular, complexes of chains and cochains with some further structures imposed on them, as well as the concepts of ellipticity and naturality. Moreover, we address homology and cohomology of such complexes. In section~3, we examine the correspondence between cochains and differential forms. In section~4, we discuss how to achieve a discrete formulation of a partial differential equation description of continuum physics in the context of FEM. Any background in differential geometry and algebraic topology is useful. For this, we refer the reader to e.g. \citep{Gross} and \citep{Frankel}. In section~5, we discuss natural differential operators in the discrete setting. Here, we present our main result (Theorem \ref{theorem1}). We first define discrete vector fields via a perturbation of the vertices of a simplex, and using this definition, we show that a discrete Lie derivative is obtained from an analogue of Cartan's magic formula. 
Finally, in section~6, conclusions are drawn.


\section{Mathematical structures for discrete solvers}

Physical quantities are related to geometrical objects. For example, we have force-like quantities related to paths, flows and fluxes through surfaces, and quantities contained by a volume. Formalizing this intuition in terms of mathematics, these geometrical objects are chains, formal linear combinations of what we call cells, and the linear functionals on such spaces are called cochains. A thorough exposition of such algebraic topological concepts in the context of electromagnetic modelling can be found in \citep{Gross}. Ker\"anen's thesis \citep{Keranen2011} offers an accessible treatment of cochains, although it has an emphasis on electromagnetism on \emph{semi}-Riemannian manifolds. For a classic exposition, we refer the reader e.g. to \citep{Whitney}.

Here, we shall recall the notions of chain and cochain complexes, and the relevant complexes with additional structures, elliptic complexes in particular. We will be rigorous but brief, leaving some details to be elaborated on in the references.

\subsection{(Co)chain complexes.}

To be able to model quantities related to geometric objects of different dimensions, we want to combine the (co)chains of different dimension in a single structure, a (co)chain complex on a topological space $\Omega$, a sequence of Abelian groups or modules. As is often desirable in continuum physics, dealing with (co)chains with coefficients in $\mathbb{R}$, the resulting spaces are vector spaces. 

Combining the information we have in the spaces of chains of different dimensions $C_i(\Omega)$ and the boundary operators $\partial$ acting as mappings between them, we obtain the sequence
\begin{equation}\label{chainComplex}
	0 \rightarrow C_n(\Omega) \stackrel{\partial}{\rightarrow} C_{n-1}(\Omega) \stackrel{\partial}{\rightarrow}... \stackrel{\partial}{\rightarrow} C_0(\Omega) \rightarrow 0,
\end{equation}
with $\partial \circ \partial = 0$, defining a {\bf chain complex} on $n$-dimensional $\Omega$, $C_*(\Omega)$. Dually, a {\bf cochain complex} on $\Omega$, $C^*(\Omega)$, is the sequence
\begin{equation}\label{cochainComplex}
	0 \leftarrow C^n(\Omega) \stackrel{\rmd}{\leftarrow} C^{n-1}(\Omega) \stackrel{\rmd}{\leftarrow}... \stackrel{\rmd}{\leftarrow} C^0(\Omega) \leftarrow 0
\end{equation}
of cochain spaces of different dimensions $C^i(\Omega)$ with coboundary operators $\rmd$ between them. Again, $\rmd \circ \rmd = 0$. Such complexes are {\bf exact} if the images of the (co)boundaries match with the kernels of the following ones. The familiar Stokes' theorem \cite{Frankel} is a nondegenerate bilinear pairing between chains and cochains.

A {\bf subcomplex of a chain complex}, or a {\bf chain subcomplex} of $C_*(\Omega)$ is a chain complex $\kappa_*(\Omega)$ consisting of spaces $\kappa_p(\Omega) \subset C_p(\Omega)$ with  $\partial \kappa_{p+1}(\Omega) \subset \kappa_{p}(\Omega)$ \citep{Naber}. Similarly, a {\bf cochain subcomplex} $\kappa^*(\Omega)$ of $C^*(\Omega)$ consists of spaces $\kappa^p(\Omega) \subset C^p(\Omega)$ with $\rmd \kappa^p(\Omega) \subset \kappa^{p+1}(\Omega)$. Note that a finite (co)chain complex arising from a cellular mesh complex on $\Omega$ is a subcomplex of the infinite-dimensional (co)chain complex on $\Omega$.

A {\bf chain map} between chain complexes $C_*(\Omega_1)$ and $C_*(\Omega_2)$ is a sequence of structure-preserving mappings $f_i: C_{i}(\Omega_1) \rightarrow C_{i}(\Omega_2)$ for all $i \in \{0,1,...,n\}$ such that $f_{i} \circ \partial = \partial \circ f_{i+1}$. That is, the boundary maps on the two complexes commute with $f_i$. Dually a {\bf cochain map} between cochain complexes $C^*(\Omega_1)$ and $C^*(\Omega_2)$ is a sequence of structure-preserving mappings $f^i: C^{i}(\Omega_1) \rightarrow C^{i}(\Omega_2)$ for all $i \in \{0,1,...,n\}$ such that $f^{i} \circ \rmd = \rmd \circ f^{i-1}$. (Co)chain maps are thus induced by mappings between $\Omega_1$ and $\Omega_2$. (Co)chain complexes and (co)chain maps form the {\bf category of (co)chain complexes}. $C_*(\Omega)$ and $C^*(\Omega)$ are functors from the category of topological spaces and continuous maps to this category, attaching a (co)chain complex to each topological space and a (co)chain map to each continuous mapping between such spaces. For relevant definitions of categories, functors, and natural transormations for the purposes of this paper, we refer the reader to \citep{Adamek}.

\subsection{Elliptic complexes}

The \emph{elliptic complex} was popularized over half a century ago by Donald Spencer, primarily through his student J. J. Kohn \citep{Kohn}. Spencer's 1969 paper ties elliptic complexes to homological algebra and overdetermined PDEs \citep{Spencer69}. Advantageously, an elliptic complex gives rise to a Hodge theory which, when restricted to the de Rham complex, yields the traditional Hodge theory of differential forms (see e.g. \citep{Wells} for a more recent exposition).  

To articulate what exactly is an elliptic complex, we need a few definitions. We shall follow the notation of \citep{Gilkey}, which is a standard reference for more details. Confining ourselves first to $\mathbb{R}^n$, with $x = (x_1,...,x_n) \in \mathbb{R}^n$, and all functions complex-valued, a {\bf linear differential operator of order $m$} $P(x,D) = \sum_{{|\alpha|} \leq m} a_\alpha (x) D^{\alpha}_{x}$ has the associated {\bf symbol} $\sigma(P) = p(x, \chi) = \sum_{|\alpha| \leq m} a_\alpha(x)\chi^\alpha$, with formal variables $\chi^\alpha$ replacing the partial derivatives $D^{\alpha}_{x}= (-\mathrm{i})^{|\alpha|}{\left(\frac{\partial}{\partial x_1}\right)}^{\alpha_1} ... \quad {\left(\frac{\partial}{\partial x_n}\right)}^{\alpha_n}$ with respect to the multi-index $\alpha$, an ordered tuple of indices. Here, $|\alpha| = \sum_{i=1}^n{ \alpha_i}$, and $\mathrm{i}$ is the imaginary unit. Moreover, interpreting the pair $(x, \chi)$ as an element of the cotangent bundle, essentially, the symbol replaces partial derivatives with coordinates of the cotangent bundle. The highest-degree part of the symbol is called the {\bf leading-order symbol} $\sigma_\mathrm{L}P = \sum_{|\alpha| = m} a_\alpha(x)\chi^\alpha$. 

By inverse Fourier transform and the properties of the Fourier transform we can write 
\begin{equation}\label{PDO}
	P(x,D)u(x) = \int_{\mathbb{R}^n}{e^{ix \chi} p(x,\chi)\hat{u}(\chi) \rmd \chi},
\end{equation}
where $\hat{u}$ denotes the Fourier transform of $u$. Then, a {\bf pseudo-differential operator} $P(x,D)$ is defined to operate on the function $u(x)$ by \eqref{PDO}. Here, $p(x,\chi)$, a smooth function in $x$ and $\chi$ with compact $x$-support, is a {\bf symbol of order $m$}, such that for all multi-indices $\alpha$ and $\beta$ there exists a constant $C_{\alpha \beta}$, such that $|D^{\alpha}_{x}  D^{\beta}_{\chi}p(x,\chi)| \leq C_{\alpha \beta}(1+|\chi|)^{m-|\beta|}$. Then, the corresponding pseudo-differential operator is said to be a {\bf pseudo-differential operator of order $m$}. Note how, essentially, \eqref{PDO} is composed of Fourier transform, multiplication by the symbol function and inverse Fourier transform. Thus, one may note the relation to \emph{filtering} in signal processing. 

Pseudo-differential operators form an algebra closed under composition and adjoints. An {\it elliptic operator} is essentially an invertible element in this algebra. More precisely a pseudo-differential operator is an {\bf elliptic operator} if its symbol $p(x,\chi)$ is elliptic, which means that $p$ is invertible and there exists $C$ such that $|p(x,\chi)^{-1}| \leq C(1+|\chi|)^{-m}$, for a $\chi$ large enough.

Finally, we note that all we have said above about pseudo-differential and elliptic operators generalizes nicely to a compact Riemannian manifold $\Omega$ by working with local coordinate charts \citep{Gilkey}. That is, an operator on smooth functions on a compact Riemannian manifold is a pseudo-differential operator, if the induced map on coordinate charts is a pseudo-differential operator for each chart $\psi_\alpha: U_\alpha \subset \Omega \rightarrow \mathbb{R}^n$. Moreover, this generalizes to vector bundles\footnote{Why consider vector bundles in the first place? If we are to discuss field quantities, we are to discuss vector bundles. The field quantities we are solving for when modelling engineering problems are sections of vector bundles. For example, a vector field is a section of the tangent bundle, and a covector field, a differential 1-form, is a section of the cotangent bundle.} over $\Omega$: An operator $P: \Gamma(E_1) \rightarrow \Gamma(E_2)$ is a pseudo-differential operator of order $m$ if $P$ is locally expressible as a matrix, each component of which is a pseudo-differential operator of order $m$. Here, $\Gamma$ is the functor taking a vector bundle $E_i$ to the space of its smooth sections. Then, for each covector $\chi \in \TstarxM$, $\sigma(P): \TstarxM \rightarrow \mathrm{Hom}(E^{x}_{1}, E^{x}_{2})$ can be written as a matrix. The symbol of $P$ can thus be viewed as a mapping taking covectors to morphisms between the fibers of $E_1$ and $E_2$ over $x$. The operator $P$ being elliptic translates to the leading order symbol $\sigma_\mathrm{L}P$ mapping $E^{x}_{1}$ isomorphically to $E^{x}_{2}$, for all $x \in \Omega$ and $0 \neq \chi \in \TstarxM$.  

Then, consider the sequence of spaces of sections of vector bundles \cite{Terng} $E_i$ over $\Omega$ and differential operators $D$ locally expressible as matrices of partial derivatives:
\begin{equation}\label{EllipticComplex}
	\Gamma(E_n) \stackrel{D}{\leftarrow} ... \stackrel{D}{\leftarrow} \Gamma(E_1) \stackrel{D}{\leftarrow} \Gamma(E_0).
\end{equation}
At each point $x$, we get the associated \emph{sequence of symbols}
\begin{equation}\label{EllipticSymbolSequence}
	0 \leftarrow E^{x}_{n}  \stackrel{\sigma_\mathrm{L}D}{\leftarrow} ... \stackrel{\sigma_\mathrm{L}D}{\leftarrow} E^{x}_{1} \stackrel{\sigma_\mathrm{L}D}{\leftarrow} E^{x}_{0} \leftarrow 0.
\end{equation}

The sequence \eqref{EllipticComplex} is a complex if $D \circ D = 0$ and it is an {\bf elliptic complex} if \eqref{EllipticSymbolSequence} is exact for each $x$ and for each $\chi \neq 0$. Moreover, for boundary conditions yielding a well-posed boundary value problem (the Shapiro-Lopatinskii boundary conditions \citep{Shapiro}, \citep{Lopatinskii}), an elliptic complex becomes a {\bf Fredholm complex} \citep{Segal}. Fredholm operators are characterized by finite-dimensional kernel and cokernel -- the Hilbert spaces are separable, yielding a countable basis. Note that even though the Fredholm operators may act on infinite-dimensional spaces, we obtain a very tangible picture of such a complex in terms of the symbol sequence via finite-dimensional vector spaces and matrices. Computations are formulated in terms of finite-dimensional complexes coming from the symbol sequence. Moreover, all this should be rather intuitive for an engineer by a ``generalization of the Laplace transform''.

Many elliptic operators of interest arise from elliptic complexes. From ellipticity of the complex, it follows that the sum of the differential operator and its adjoint $D + D^*$ is an elliptic operator \citep{Atiyah}. Given an elliptic complex, the elliptic operator $\Delta = DD^* + D^*D$ yields the harmonic sections of the complex, analogously to harmonic forms of the de Rham complex, and a generalized Hodge theory in this manner. Thus, with an elliptic complex, we are provided with an extremely useful set of tools for doing computational physics.

In general, complexes of {\bf Hilbert spaces}, i.e. vector spaces equipped with an inner product that are complete with respect to the induced norm have been recently in the forefront within the context of numerical methods \citep{Arnold, ArnoldBook, Hiptmair2020}. Such general complexes are not necessarily Fredholm or elliptic, and the spaces are not necessarily separable. A {\bf Hilbert complex} \citep{Arnold, ArnoldBook, Hiptmair2020} $\mathcal{H}^*$ is the sequence of Hilbert spaces
\begin{equation}\label{HilbertComplex}
	0 \leftarrow \mathcal{H}^n \stackrel{\rmd}{\leftarrow} \mathcal{H}^{n-1} \stackrel{\rmd}{\leftarrow}... \stackrel{\rmd}{\leftarrow} \mathcal{H}^0 \leftarrow 0,
\end{equation}
with $\rmd$ representing here closed densely-defined linear operators between Hilbert spaces, with $\rmd \circ \rmd = 0$ \citep{Arnold}, \citep{Bruning}, \citep{HolstAndStern}. An elliptic complex is a Hilbert complex. The most obvious example of a Hilbert complex of cochains is the {\bf de Rham complex} $F^*(\Omega)$ of differential forms on a smooth manifold $\Omega$, 
\begin{equation}\label{deRhamComplex0}
	0 \leftarrow F^n(\Omega) \stackrel{\rmd}{\leftarrow} F^{n-1}(\Omega) \stackrel{{\rmd}}{\leftarrow}... \stackrel{{\rmd}}{\leftarrow} F^0(\Omega) \leftarrow 0,
\end{equation}
with the exterior derivative acting as the coboundary mapping. Essentially, Hilbert complex can be viewed as an abstraction of the de Rham complex. A {\bf Hilbert subcomplex} $\mathcal{H}^{*}_\mathrm{s}$ of $\mathcal{H}^{*}$ consists of spaces $\mathcal{H}^{p}_\mathrm{s} \subset \mathcal{H}^p$ with $\rmd \mathcal{H}^{p}_\mathrm{s} \subset \mathcal{H}^{p+1}_\mathrm{s}$. For Hilbert subcomplexes, the subspaces are naturally understood as Hilbert subspaces.

However, for computational continuum physics when compact domains are considered, a Hilbert complex is often too abstract of a notion, since the more pedestrian approach of elliptic and Fredholm complexes would suffice, as we will see in the following sections. An elliptic complex, as a special case of a Hilbert complex, yields thus in many cases a more concrete approach to our needs than a general Hilbert complex. The de Rham complex is indeed an elliptic complex. As also noted in \citep{Gross} and \citep{Kotiuga1984}, it is the framework of elliptic complexes that formalizes the ideas of analogies in physics expressed by Tonti via Tonti diagrams \citep{Tonti1972}, \citep{Tonti1977}: the de Rham isomorphism at one's disposal in the complexes associated with the exterior derivative gives a concrete tool to answer questions about (co)homology of the complex.

\subsection{Naturality}

The concept of naturality is ubiquitous in mathematics. In particular, we are here interested in natural vector bundles and natural differential operators. Seminal early work in a context relevant for us includes that of Weil \citep{Weil1953} and Palais \citep{Palais}. The book by Kolar~{\it et~al.} \citep{Kolar} is a modern treatment fully set in the language of category theory. 

\subsubsection{Naturality associated with differentiable structure}

A {\bf natural vector bundle} is a functor $E$ from the category of $n$-dimensional smooth manifolds and local diffeomorphisms to the category of vector bundles and vector bundle homomorphisms. Let us denote the canonical projections of the vector bundles $E(\Omega)$ and $E(\Omega')$ as $\pi_{\Omega}$ and $\pi_{\Omega'}$. The functor $E$ associates a vector bundle to each $n$-manifold $\Omega$ and to each $f: \Omega \rightarrow \Omega'$ a vector bundle homomorphism, i.e. a commutative diagram
\[\begindc{\commdiag}[50]
\obj(0,10)[a]{$E(\Omega)$}
\obj(20,10)[c]{$E(\Omega')$}
\obj(20,0)[d]{$\Omega'$}
\obj(0,0)[b]{$\Omega$}
\mor{a}{c}{$E(f)$}[\atleft, \solidarrow]
\mor{c}{d}{$\pi_{\Omega'}$}[\atleft, \solidarrow]
\mor{a}{b}{$\pi_{\Omega}$}[\atright, \solidarrow]
\mor{b}{d}{$f$}[\atleft, \solidarrow]
\enddc\]
which is a fiberwise linear isomorphism and covers $f$ \citep{Kolar}. Many typical vector bundle constructions, such as bundles of tangent $p$-vectors and $p$-covectors, are natural: morphisms of certain type between manifolds are lifted to morphisms between vector bundles over the manifolds. For example, in the case of tangent bundles, in fact any smooth function between manifolds may be lifted to its push-forward between the corresponding tangent bundles. In this case, interpreting tangent vectors as differential operators, functoriality manifests itself in the chain rule. Moreover, a {\bf natural differential operator} between natural vector bundles is, essentially, such that it commutes with local diffeomorphisms \citep{Stredder}, \citep{Terng}, \citep{Kolar}. Moreover, if the associated vector bundles of an elliptic complex are \emph{natural vector bundles} and the operators \emph{natural differential operators}, we call the elliptic complex a {\bf natural elliptic complex}.

For example, the exterior derivative $\rmd$ of differential forms is a natural differential operator as it commutes with pullbacks by smooth functions. More precisely, let $\mathcal{D}F^{p}(\Omega)$ the space of smooth $p$-forms on $\Omega$ i.e. the space of smooth sections of the (natural) vector bundle of $p$-covectors, and $f^*$ the pullback of a smooth map $f$ between the smooth manifolds $\Omega'$ and $\Omega$. Then, the square
\[\begindc{\commdiag}[50]
\obj(0,10)[a]{$\mathcal{D}F^{p}(\Omega)$}
\obj(20,10)[c]{$\mathcal{D}F^{p}(\Omega')$}
\obj(20,0)[d]{$\mathcal{D}F^{p+1}(\Omega')$}
\obj(0,0)[b]{$\mathcal{D}F^{p+1}(\Omega)$}
\mor{a}{c}{$f^*$}[\atleft, \solidarrow]
\mor{c}{d}{$\rmd$}[\atleft, \solidarrow]
\mor{a}{b}{$\rmd$}[\atright, \solidarrow]
\mor{b}{d}{$f^*$}[\atleft, \solidarrow]
\enddc\]
commutes. This is to say that the exterior derivative is a natural transformation between the (contravariant) functors that assign the spaces of smooth $p$-forms and $(p+1)$-forms to $\Omega$, respectively. In fact, it is essentially \emph{the unique natural differential operator} between the associated natural vector bundles \citep{Palais}, \citep{Kolar}. As pointed out in \citep{Kolar}, linearity of $\rmd$ is thus also a consequence of naturality.

Another important natural operation commuting with pullbacks is the contraction of a differential $p$-form $\eta$ with a vector field $\bf{v}$, $\mathrm{i_{\bf v} \eta}$, producing a $(p-1)$-form. As an example, the contraction of magnetic flux density with the velocity field {\bf v}, $\mathrm{i_{\bf v} B}$ is a natural way to interpret the ${\bf v} \times {\bf B}$ term of the Lorentz force, essential in, e.g., modelling magnetohydrodynamics (MHD). Through Cartan's magic formula 
\begin{equation}\label{cartan}
	\mathcal{L}_{\bf v} = \mathrm{i_{\bf v}} \circ \rmd + \rmd \circ \mathrm{i_{\bf v}},
\end{equation}
a structure-preserving discretization of contraction is also important for finding a discrete analogue for the Lie derivative $\mathcal{L}_{\bf v}$. Bossavit's take on discretizing contractions is via the duality between the extrusion of a manifold (extruding along the flow of a vector field) and contraction of a differential form, utilizing Whitney forms \citep{Bossavit03}. In section~5, we will give our discrete version of Cartan's magic formula yielding a discretization for the Lie derivative (Theorem 1).

\subsubsection{Naturality associated with geometric structure}

Different type of natural operations arise when a geometric structure, such as an inner product, on a vector bundle is considered \citep{Kolar}. Given an inner product, we have the associated connection and curvature, and the associated natural operations involve preservation of \emph{isometries} and not just diffeomorphisms. In a sense, we need to restrict the category we are working in. We will see manifestations of such naturality later on in the context of FEM in section~4, manifest in the Galerkin step.

Such naturality is not, however, restricted to the Galerkin step. If we wanted to discuss mixed tensor physics, we would necessarily need more complicated operators than just the exterior derivative. As an example, small-strain elasticity may be articulated using the covariant exterior derivative $\rmd_\nabla$ and vector- and covector-valued differential forms as discussed recently e.g. in \citep{Kovanen}.\footnote{For a thorough exposition on the theory of elasticity, we refer to \citep{Marsden}.} In general $\rmd_\nabla \circ \rmd_\nabla \neq 0$, unless one is dealing with a flat connection. Thus, when seeking naturality related to geometric structure, the notion of curvature is immediate. 

Concretely, the connection compatible with the Riemannian metric, the Levi-Civita connection $\nabla$ is natural with respect to isometries. In terms of vector fields, this means that given $\hat{f}: \Omega \rightarrow \Omega'$ an isometry between Riemannian manifolds with Levi-Civita connections $\nabla$ and $\nabla'$, respectively, and ${\bf v} \in \TxM$ a tangent vector at a point $x$ and ${\bf w} \in \Gamma(\TM)$ a vector field i.e. a section of the tangent bundle $\TM$, $\hat{f}_*(\nabla_{{\bf v}} {\bf w}) = \nabla'_{\hat{f}_*{\bf v}} (\hat{f}_* {\bf w})$, where $\hat{f}_*$ denotes the pushforward along $\hat{f}$, holds. Diagrammatically, the square
\[\begindc{\commdiag}[50]
\obj(0,10)[a]{$\Gamma(\mathrm{T}\Omega)$}
\obj(20,10)[c]{$\Gamma(\mathrm{T}\Omega')$}
\obj(20,0)[d]{$\Gamma(\mathrm{T}\Omega')$}
\obj(0,0)[b]{$\Gamma(\mathrm{T}\Omega)$}
\mor{a}{c}{$\hat{f}_*$}[\atleft, \solidarrow]
\mor{c}{d}{$\nabla_{\hat{f}_*{\bf v}}'$}[\atleft, \solidarrow]
\mor{a}{b}{$\nabla_{\bf v}$}[\atright, \solidarrow]
\mor{b}{d}{$\hat{f}_*$}[\atleft, \solidarrow]
\enddc\]
commutes. For proof, see e.g. \citep{Poor}. Considering structure-preserving (FEM) discretizations of e.g. elasticity and fluid mechanics, the naturality of the connection, or the covariant derivative, to be preserved thus involves invariance with respect to the Riemannian structure \emph{a priori}, before the Galerkin step, unlike e.g. in the case of discrete electromagnetism. To tackle such physics and its natural operations in a structure-preserving manner in the discrete setting, one needs to take into account the intimate tie between the natural differential operators used and the Riemannian structure. A systematic approach acknowledging the structures behind the natural operations, keeping in mind the challenges of a DEC or FEEC type of discretization \citep{Hirani}, \citep{Arnold}, is bound to face constraints \citep{Kotiuga2008}.

We conclude this subsection with two open research problems.

\begin{myrp}[Discrete natural operations in elasticity and fluids]
	
	The Euler equations of fluid dynamics involve the naturality of the de Rham complex, while elasticity and the viscosity term of the Navier-Stokes equations involve invariance under structure-preserving isometries. Naturality unifies and ties categorical notions to discretizations of such equations via Andre Weil's 1953 paper \citep{Weil1953} and the modern setting exemplified by \citep{Kolar}, which captures the seemingly different aspects of naturality under the same structure. Keeping in mind that the discrete exterior calculus \citep{Hirani} and finite element exterior calculus \citep{Arnold} have their limitations in the context of the commutative cochain problem and the discrete star localization problem \citep{Kotiuga2008}, a systematic approach to a structure-preserving discretization needs to take into account the inherent difference in the type of naturality of differential operators.
\end{myrp}

\begin{myrp}[Natural differential operators in multiphysics]
	
	The different notions of naturality related to e.g. MHD, plasmas and superconductivity should be studied. These involve geometric free boundary problems associated with them. These in turn relate to near force-free magnetic fields, having implications for e.g. fusion research.
	
\end{myrp}


\subsection{Homology and cohomology}

For modelling needs, chain and cochain complexes provide us with a way to discuss topological invariants. In particular, \emph{homology} and \emph{cohomology} capture the quantity and quality of holes in $\Omega$. Intuitively speaking, take a 3-dimensional domain $\Omega$ for example. There, 0-homology captures the quantity of connected components of $\Omega$, 1-homology the number of tunnels through $\Omega$ and 2-homology the number of voids. Cohomology can then be viewed as assigning values to these holes. In CEM, (co)homological considerations of the modelling domain $\Omega$ are at the very heart of the science of modelling \citep{Kotiuga1984, Gross}; Consider for example, driving net currents through tunnels in the modelling domain in a magneto(quasi)static problem.\footnote{Exploiting cohomology can be a significant difference-maker in terms of efficiency of computations as well. As an example, see e.g. \citep{LahtinenJSNM} for making use of cohomology in non-linear magnetoquasistatic problems of superconductor modelling.}

In terms of the chain complex and the kernel and the image of the boundary mapping $\partial$, the {\bf $p$th homology group} may be defined as
\begin{equation}\label{homology}
	H_p = \mathrm{ker} (\partial_p) / \mathrm{im} (\partial_{p+1}),
\end{equation}
i.e. the quotient group formed by $p$-dimensional {\bf cycles} (chains with empty boundaries) in $C_*(\Omega)$ modulo $(p+1)$-dimensional {\bf boundaries} ($(p+1)$-chains that are boundaries for a $p$-chain) in $C_*(\Omega)$. Similarly, we can define the {\bf $p$th cohomology group} as
\begin{equation}\label{cohomology}
	H^p = \mathrm{ker} (\rmd_p) / \mathrm{im} (\rmd_{p-1}),
\end{equation}
i.e. the quotient group formed by $p$-dimensional {\bf cocycles} (cochains whose coboundary vanish) in $C^*(\Omega)$ modulo $(p-1)$-dimensional {\bf coboundaries} ($(p-1)$-cochains that are coboundaries for a $p$-cochain) in $C^*(\Omega)$. Thus, a complex is exact if and only if its (co)homology vanishes.

(Co)homology theories with seemingly different foundations abound. For example, cellular homology, simplicial homology, Cech homology and de Rham homology are homology theories dubbed under different labels, to name just a few. However, since the Eilenberg-Steenrod\footnote{Relaxing one of these axioms, the \emph{dimension axiom}, yields a {\bf generalized (co)homology theory}. Two standard examples are called topological $K$-theory and algebraic \emph{$K$-theory} \citep{Swan}.} axioms were laid out, they can all be seen as instances of the general theory, and it turns out that homology is a homotopy invariant of topological spaces \citep{Eilenberg-Steenrod}. The Eilenberg-Steenrod axioms, of which there are versions for both homology and cohomology (obtained from one another by reversing the arrows), assert that a homology theory is a sequence of functors from a suitable category of topological spaces (or from the category of chain complexes) to category of Abelian groups (or sequences of them), together with boundary operators $\partial$, which induce \emph{natural transformations between homology functors}. This system is required to satisfy five properties, for which the reader is referred to \citep{Eilenberg-Steenrod}. These axioms allow us to talk about (co)homology theories without referencing to any concrete (co)chain realization.

But how do we view $H_p$ as functors? First of all, the assignment of a chain complex to the space $\Omega$ is functorial. Then, the {\bf homology} is a (covariant) functor from the category of chain complexes to the category of sequences of Abelian groups and sequences of Abelian group homomorphisms, satisfying the Eilenberg-Steenrod axioms. Dually, {\bf cohomology} is a contravariant functor of the same kind from the category of cochain complexes. For homology, the functoriality means here, that given two chain complexes and a chain map $f$ between them, homology functor takes the chain spaces to the homology groups and the chain map to induced homomorphisms between the groups (as chain maps take cycles to cycles and boundaries to boundaries), in a way that respects identity maps and composition. Then, there is a composite functor which first assigns a chain complex to a space and then the homology groups to the complex. This is the homology of $\Omega$.

Moreover, the assignment of the $p$th homology group to a space may be viewed as a functor in its own right. Then, the boundary operators $\partial$ induce natural transformations between homology functors in the sense that for all $0<p \leq n$ they assign to each $\Omega$ a homomorphism $\partial: H_p(\Omega) \rightarrow H_{p-1}(\Omega)$ such that for each morphism $f$ from $\Omega$ to $\Omega'$ the diagram
\[\begindc{\commdiag}[50]
\obj(0,10)[a]{$H_p(\Omega)$}
\obj(0,0)[c]{$H_{p-1}(\Omega)$}
\obj(20,0)[d]{$H_{p-1}(\Omega')$}
\obj(20,10)[b]{$H_p(\Omega')$}
\mor{a}{c}{$\partial$}[\atright, \solidarrow]
\mor{c}{d}{$H_{p-1}(f)$}[\atright, \solidarrow]
\mor{a}{b}{$H_{p}(f)$}[\atleft, \solidarrow]
\mor{b}{d}{$\partial'$}[\atleft, \solidarrow]
\enddc\]
\noindent commutes.\footnote{Notation is slightly abused here by denoting the induced boundary homomorphisms with the same symbol as the boundary operators of chain complexes.} The existence of such a natural transformation can be articulated via the zig-zag (or snake) lemma which yields a commutative diagram of long exact sequences of homology groups and group homomorphisms from a commutative diagram of short exact sequences of chain complexes and chain maps (see e.g. \citep{Ghrist}). Note that even though we have considered the homology of chain complexes here for concreteness, this functoriality is at the heart of every homology theory. The dual results hold for cohomology: then, the coboundary operators $\rmd$ induce such natural transformations.

As a prelude to section~3, we point out that for the de Rham complex, we have the \emph{de Rham cohomology}, where the role of cocycles is played by closed differential forms (those whose exterior derivative vanish) and exact forms (those that are exterior derivatives of another form) are the coboundaries. Then, utilizing the induced Hodge theory, the representatives of cohomology classes are harmonic forms: There is a unique harmonic differential $p$-form with prescribed periods on the homology basis for $p$-cycles on $\Omega$ with coefficients in $\mathbb{R}$. As the de Rham isomorphism can be utilized, cohomology of this complex becomes tangible, and Tonti diagrams obtain a clear formal meaning \citep{Gross}. Clearly, the same thing cannot be said about a general Hilbert complex \citep{Arnold}. However, for a general elliptic complex, it is guaranteed that we have an analogous cohomology theory: The space of harmonic sections, the kernel of $\Delta = DD^* + D^*D$, of the complex is isomorphic to the cohomology of the complex of corresponding degree \citep{Gilkey}. That is, the dimension of the space of harmonic sections of order $p$ on $\Omega$ is the $p$th Betti number of $\Omega$. Thus, elliptic complexes are the framework to 	formalize Tonti's analogies in physics.

\subsection{An electrical engineering interlude}

The relevant data structure for computational physics of continuum is that of a (co)chain complex. With the additional structure of an elliptic complex and naturality of the associated differential operators, one has a machine that produces a Hodge theory, extremely effective for articulating physical theories and essential for numerical methods in continuum physics. This framework is also essential in formalizing Tonti diagrams \citep{Gross}. Moreover, (co)homological information of the complex as a topological invariant is essential for computations. Thinking in terms of traditional electrical engineering, via Hodge theory, homology is the bridge between field theory and circuit theory: The periods of harmonic forms are the variables that appear in Kirchoff's laws, while the lumped parameters have a Hilbert space interpretation \citep{Kotiuga1984}. Thus, elliptic and Fredholm complexes, naturality and (co)homology are at the very heart of engineering. In this light, a general Hilbert complex often appears as an excessive abstraction in compact domains, without the benefits of ellipticity or separability \citep{Arnold}. Note also how the concept of naturality comes up at different levels. Desirably, such structures should be preserved in discretizations, too.

\section{Correspondence between cochains and differential forms}

The basic question addressed by Andr\'e Weil in 1952 was the fact that cohomology theories satisfying the Eilenberg-Steenrod axioms \citep{Eilenberg-Steenrod} had isomorphic cohomology groups \citep{Weil}. However, the axiomatic method made no mention of chains or cochains, and in the process, had nothing to say about possible correspondences between cochains in concrete realizations of cohomology theories such as those of de Rham and Cech, or simplicial cohomology theory. In the process of devising tools for addressing this issue, Weil developed the spectral sequence and an interpolation formula now known as Whitney forms. Hassler Whitney developed these tools for the purposes of analysis on manifolds, and published his results in his monograph entitled Geometric Integration Theory in 1957 \citep{Whitney}.

Let us scrutinize the following diagram and concentrate on the simplicial case:
\[\begindc{\commdiag}[50]
\obj(0,20)[a]{Cochains}
\obj(30,20)[b]{Differential forms}
\obj(2,21)[c]{}
\obj(23,21)[d]{}
\mor{a}{b}{Limit process}[\atright, \solidarrow]
\mor{d}{c}{Integration}[\atright, \solidarrow]
\enddc\]
Given a simplicial complex $K$ on a smooth manifold $\Omega$, the \emph{Whitney map} $W$ interpolates a cochain in a piecewise linear manner creating a Whitney form. Dually, the \emph{de Rham map} $R$ integrates a differential form $\eta$ to provide a cochain. The composite $RW$ is the identity mapping on cochains. Moreover, there is a precise sense in which $\eta$ is homotopic to $WR\eta$, and the composite $WR$ tends to identity with the refinement of the simplicial complex \citep{Dodziuk1976}.

Formally, given the barycentric coordinates $\mu_{i}$ of the vertices in a simplicial complex $K$, the {\bf Whitney map} $W$, which takes a simplicial $p$-cochain to a differential $p$-form, gives a Whitney $p$-form on a $p$-simplex $\sigma = \{i_0, i_1, ... i_p\}$ via
\begin{equation}\label{WhitneyMap}
	W \sigma = p!\sum_{k=0}^p(-1)^k{\mu_{i_k}}\rmd \mu_{i_0} \wedge ... \wedge \widehat{\rmd \mu_{i_k}} \wedge ... \wedge \rmd \mu_{i_p}, \quad p > 0.
\end{equation}
For $p=0$ we have simply $W\sigma = \mu_{i_0}$. Here, $\wedge$ is the exterior product of differential forms, which, given a $p$-form and a $q$-form, produces a $p+q$-form, and $\widehat{\cdot}$ denotes that we omit the factor under the hat sign. Moreover, the {\bf de Rham map} $R$ takes a differential $p$-form $\eta$ to a simplicial $p$-cochain by integrating it on a $p$-chain $c$:
\begin{equation}\label{deRhamMap}
	R \eta = \int_c \eta.
\end{equation}

Functors $C^*$ and $F^*$ assign the {\bf simplicial cochain complex} $C^*(K,\Omega)$ and the de Rham complex $F^*(\Omega)$ to $\Omega$. Then, the Whitney map is a cochain map between the simplicial cochain complex
\begin{equation}\label{SimplicialComplex}
	0 \leftarrow C^n(K, \Omega) \stackrel{\rmd}{\leftarrow} C^{n-1}(K, \Omega) \stackrel{\rmd}{\leftarrow}... \stackrel{\rmd}{\leftarrow} C^0(K, \Omega) \leftarrow 0
\end{equation}
and the de Rham complex
\begin{equation}\label{deRhamComplex1}
	0 \leftarrow F^n(\Omega) \stackrel{\rmd}{\leftarrow} F^{n-1}(\Omega) \stackrel{{\rmd}}{\leftarrow}... \stackrel{{\rmd}}{\leftarrow} F^0(\Omega) \leftarrow 0.
\end{equation}
Moreover, the de Rham map is a cochain map in the other direction. Now, $W \circ \rmd \tilde{c} = {\rmd} \circ W \tilde{c}$, where $\tilde{c}$ is a simplicial cochain and $\rmd$ on the left (and in \eqref{SimplicialComplex}) denotes the coboundary on simplicial cochains (the transpose of the simplicial boundary operator), while $\rmd$ on the right (and in \eqref{deRhamComplex1}) denotes the exterior derivative of a differential form.\footnote{We will make this definition natural and more explicit in section~5.} Thus, the isomorphism of de Rham cohomology groups of $\Omega$ and simplicial cohomology groups related to $K$ is implied, as cochain maps induce maps on cohomology functorially \citep{Gross}, \citep{Dodziuk1976}, \citep{Muller}. Note that the coboundary mappings then induce natural transformations on cohomology. 

All this suggests an interpretation for Whitney forms and the correspondence of cochains and differential forms as parts of a higher-categorical structure. The fact that \emph{chain homotopic} chain maps induce the same maps on homology is the core of the isomorphism of de Rham and simplicial cohomology groups \citep{BottTu}. In fact, one could say that such chain homotopy was the reason Whitney forms were conceived \citep{Weil}. By considering chain complexes, chain maps and chain homotopies between them, one obtains a 2-category with chain homotopies acting as morphisms between morphisms: 2-morphisms. Furthermore, one can consider chain homotopies between chain homotopies, iterating this construction to infinity to obtain an $\infty$-category with an infinite number of ``layers'' of morphisms.\footnote{Such considerations lead us to the foundations of \emph{homotopy type theory} \citep{HoTT}. This points to a plethora of possibilities of such tools in the simplicial context \citep{Getzler}.}

We conclude this section with an open research problem.

\begin{myrp}[Whitney forms in a general category theoretical framework]
	
	The Eilenberg-Steenrod axioms yield means of dealing with essential topological aspects of physics on a level separate from any cochain realization \citep{Eilenberg-Steenrod}. Andr\'e Weil's reaction to this led to the birth of Whitney forms \citep{Weil}. Given the connection of Whitney forms and cochains to rational homotopy theory \citep{SullivanRHT}, \citep{Wilson}, the axiomatics and the importance of chain homotopy in this context suggest a niche for a categorical interpretation. In particular, to provide a framework for the Whitney and de Rham maps, one should consider e.g. functors between exterior algebra bundles, cochain complexes and manifolds. 
	
\end{myrp}

\section{The path to discrete: Finite element method}

Now, we are ready to discuss how to achieve a discrete formulation of physics. We will use the familiar FEM as an example. In this section, we will introduce the relevant details of FEM.

The Galerkin FEM (from now on we will refer to Galerkin FEM as just FEM) is likely the most utilized numerical solution method of partial differential equations in natural sciences. FEM is a discretization method: A way to discretize the function space from which the solution is sought. This is achieved by finding a finite-dimensional Hilbert subspace $\mathcal{H}_\mathrm{s}$ of the Hilbert space $\mathcal{H}$ in which the solution of the original problem resides. Or, as formalized in \citep{Arnold} and \citep{HolstAndStern}, it abstracts to finding a finite-dimensional subcomplex of a Hilbert complex on a compact manifold $\Omega$. However, in such compact domains with the Shapiro-Lopatinskii boundary conditions \citep{Shapiro}, \citep{Lopatinskii}, we are indeed dealing with Fredholm complexes and thus the theory of elliptic complexes guarantees that the Hodge theory we need falls in our laps, and we have a complex of separable spaces for the needs of functional analysis, as discussed in section~2.

\subsection{Weak formulation}

Now we consider how to obtain the relevant Hilbert spaces and the weak formulation utilized in FEM. For simplicity, we will only consider elliptic problems here. For further details, see e.g. \citep{Arnold}, \citep{Bossavit1998} and \citep{Brenner}. For a thorough introduction to the required Sobolev spaces in the context of electromagnetism see \citep{KurzAuchmann}.

In a Hilbert space $\mathcal{H}$ with inner product  $\langle \cdot, \cdot \rangle$, it holds that \citep{Yosida}
\begin{equation}\label{FEMTheorem}
	\gamma \in \mathcal{H}, \quad \gamma = 0 \quad \Leftrightarrow \quad \langle \gamma, \gamma' \rangle = 0, \quad \forall \gamma' \in \mathcal{H}.
\end{equation}
This is the backbone of FEM: using the inner product of $\mathcal{H}$ we can test whether an equation holds. In FEM, we approach a discrete formulation from the typical direction, i.e., we are discretizing a problem formulated in terms of partial differential equations. We consider the equations to be formulated in terms of differential forms. Introducing the standard $L^2$ inner product in the space of piecewise smooth differential $p$-forms ($\langle \eta, \gamma \rangle = \int_\Omega \eta \wedge \star \gamma$), we can define the $L^2$ Hilbert space of such forms $L^2F^p(\Omega)$, for whose elements $\gamma$, $\langle \gamma, \gamma \rangle < \infty$ holds. Then, we define the Sobolev spaces of differential forms
\begin{equation} \label{SobolevSpace}
	L^2F^p(\tilde{\rmd}, \Omega) = \left\{\gamma \in L^2F^p(\Omega) \quad | \quad \tilde{\rmd} \gamma \in L^2F^{p+1}(\Omega)\right\}
\end{equation}
and
\begin{equation} \label{CoSobolevSpace}
	L^2F^p(\tilde{\delta}, \Omega) = \left\{\gamma \in L^2F^p(\Omega) \quad | \quad \tilde{\delta} \gamma \in L^2F^{p-1}(\Omega)\right\},
\end{equation}
where $\tilde{\rmd}$ and $\tilde{\delta}$ denote the \emph{weak exterior derivative} and the \emph{weak co-derivative}, respectively. The Hodge operator $\star$ (see e.g. \citep{Frankel} section~14), induced by the metric on $\Omega$, is an isomorphism between $L^2F^p(\tilde{\rmd}, \Omega)$ and \linebreak $L^2F^{n-p}(\tilde{\delta}, \Omega)$. The ubiquity of $\star$ in constitutive equations of physical theories underlines the importance of both \eqref{SobolevSpace} and \eqref{CoSobolevSpace}. 

$L^2F^p(\tilde{\rmd}, \Omega)$ and $L^2F^{n-p}(\tilde{\delta}, \Omega)$ are the kind of spaces we want physical quantities to live in when using FEM: We need the quantities and their weak exterior derivatives and weak co-derivatives to be sufficiently well-behaved. The weakness of these operators ensures that they get along with piecewise smooth differential forms, too. In fact, for smooth forms, the weak and strong versions of these operators coincide. The piecewise-smoothness is often crucial in physics, as e.g. material boundaries can exhibit jumps in the derivatives of field quantities. The operation of the {\bf strong co-derivative} on a $p$-form can be defined using the Hodge operator $\star$ and the (strong) exterior derivative $\rmd$ as $(-1)^p \star^{-1} \rmd \star$. This suggests us the definition of the {\bf weak exterior derivative} $\tilde{\rmd}$, which is defined to satisfy $\langle \tilde{\rmd} \gamma, \eta \rangle =  \langle \gamma, (-1)^p \star^{-1} \rmd \star \eta \rangle, \quad \forall \eta \in \mathcal{D}F^{p+1}(\Omega)$, with vanishing boundary terms. Then, the {\bf weak co-derivative} $\tilde{\delta}$ obeys $\langle \tilde{\delta} \gamma, \varpi \rangle = \langle \gamma, \rmd \varpi \rangle, \quad \forall \varpi \in \mathcal{D}F^{p-1}(\Omega)$. Here $\mathcal{D}F^{p}(\Omega)$ is the space of smooth $p$-forms supported by $\Omega$. So, we can define the operation of these weak operators on piecewise smooth forms through operation of the strong operators on smooth forms, utilizing the inner product. Often, it is customary to make no notational difference between the weak operators and their strong counterparts, so that, for example simply $\rmd$ is used to denote both strong and weak exterior derivative.

Working in $L^2F^p(\tilde{\rmd}, \Omega)$, we can now re-state an equation of the form $\mathrm{L}\alpha = \nu$, where $\mathrm{L}$ is a linear operator, and $\alpha$ is a differential form with $\mathrm{L}\alpha \in L^2F^p(\tilde{\rmd}, \Omega)$, using \eqref{FEMTheorem}, as
\begin{equation}\label{FEMWeightedResidual}
	\langle \mathrm{L}\alpha, \gamma' \rangle = \langle \nu, \gamma' \rangle, \quad \forall \gamma' \in L^2F^p(\tilde{\rmd}, \Omega).
\end{equation}
From this {\bf weighted residual formulation} of the problem one typically approaches the \emph{weak formulation} via {\bf integration by parts}, which utilizing $\rmd$ and $\star$ can be stated as
\begin{equation}\label{integrationByParts}
	\langle \rmd \gamma, \eta \rangle = \langle \gamma, (-1)^p \star^{-1} \rmd \star \eta \rangle + \int_{\partial \Omega} \gamma \wedge \star \eta, \forall \eta \in \mathcal{D}F^{p+1}(\Omega).
\end{equation}
Thus, the metric-dependent co-derivative can be converted into the \linebreak metric-independent exterior derivative in this process. This is favourable in terms of commutative properties of the operators. Namely, the exterior derivative commutes with pullbacks by smooth functions while the co-derivative only commutes with pullbacks by isometries \citep{Gerritsma}. Commutation under pullbacks is an essential testing ground for a discretization: We want to preserve naturality. Note how the naturality of the exterior derivative is that of subsection~2.3.1 (related to differentiable structure) and the naturality of the co-derivative is that of subsection~2.3.2 (related to geometric structure).

This {\bf weak formulation} process leads to
\begin{equation}\label{FEMWeak}
	a(\alpha, \gamma') = \langle \nu, \gamma' \rangle, \quad \forall \gamma' \in L^2F^p(\tilde{\rmd}, \Omega),
\end{equation} 
where $a(\cdot, \cdot)$ is a coercive and bounded bilinear form \citep{Brenner}.

\subsection{The FEM discretization}

The solution to the weak formulation \eqref{FEMWeak} lies in an infinite-dimensional Hilbert space. The problem is yet to be discretized. As already mentioned, in FEM this is done by finding a finite-dimensional subspace for $L^2F^p(\tilde{\rmd}, \Omega)$ (supposing $\alpha$ resides there) with which to approximate the solution space. This is achieved through what is called meshing: covering of $\Omega$ with a cellular mesh complex and attachment of a set of \emph{basis functions} to the cells of desired dimension, so that the basis functions span a Hilbert subspace $W^p(\Omega)$ of $L^2F^p(\tilde{\rmd}, \Omega)$. Then, the unknown is approximated as a sum of these functions, and the weak formulation is tested with the basis of $W^p(\Omega)$, and not with all $\gamma' \in L^2F^p(\tilde{\rmd}, \Omega)$. Thus, a finite, discrete formulation of the original problem, which due to Galerkin orthogonality (see \citep[p. 58]{Brenner}), results in an optimal approximation, is obtained.

\subsubsection{Differentiable structure}

As formalized by Arnold in \citep{Arnold}, this discretization can be seen as, not only as finding a subspace of a Hilbert space, but finding a Hilbert subcomplex of the {\bf $L^2$ de Rham complex}
\begin{equation}\label{deRhamComplex}
	0 \leftarrow L^2F^n(\tilde{\rmd}, \Omega) \stackrel{\tilde{\rmd}}{\leftarrow} L^2F^{n-1}(\tilde{\rmd}, \Omega) \stackrel{\tilde{\rmd}}{\leftarrow}... \stackrel{\tilde{\rmd}}{\leftarrow} L^2F^0(\tilde{\rmd}, \Omega) \leftarrow 0.
\end{equation}
However, on a compact $\Omega$, this is indeed an elliptic complex, and the higher abstraction level of Hilbert complexes offers us no tangible benefits. Moreover, with Shapiro-Lopatinskii boundary conditions, we indeed have a Fredholm complex: a prototypical Hilbert complex with separable Hilbert spaces. Given a simplicial cellular mesh complex $K$ on $\Omega$, a viable and widely utilized option for the subcomplex is the {\bf Whitney complex}
\begin{equation}\label{WhitneyComplex}
	0 \leftarrow W^n(K, \Omega) \stackrel{{\rmd}}{\leftarrow} W^{n-1}(K, \Omega) \stackrel{{\rmd}}{\leftarrow}... \stackrel{{\rmd}}{\leftarrow} W^0(K, \Omega) \leftarrow 0,
\end{equation}
consisting of Hilbert spaces of Whitney forms, $\rmd$ here denoting the transpose of the boundary operator \citep{Bossavit}. The Whitney complex is the restriction of $L^2F^*(\tilde{\rmd}, \Omega)$ to Whitney forms, obtained via the Whitney map. At the continuum limit, as $K$ is refined, Whitney complex approaches the de Rham complex. The coboundary mappings commute with the inclusions and there is a bounded cochain mapping projecting the de Rham complex to the Whitney complex, ensuring stability. The exterior differentiation is inherited by the Whitney complex from the de Rham complex. Moreover, the naturality of the exterior derivative is preserved in this discretization: The commutative property under pullbacks still holds. 

\subsubsection{Metric-dependent properties}

In Whitney complex we have a structure-preserving discretization of the de Rham complex and the associated Hodge theory. In terms of Hodge operators, for an ideal discretization we would have a commutative diagram of the form \citep{Tarhasaari}
\[\begindc{\commdiag}[50]
\obj(0,10)[a]{$L^2F^{p}(\tilde{\rmd}, \Omega)$}
\obj(20,10)[c]{$C^p(K, \Omega)$}
\obj(20,0)[d]{$C^{n-p}(K, \Omega)$.}
\obj(0,0)[b]{$L^2F^{n-p}(\tilde{\rmd}, \Omega)$}
\mor{a}{c}{$R$}[\atleft, \solidarrow]
\mor{c}{d}{$\tilde{\star}$}[\atleft, \solidarrow]
\mor{a}{b}{$\star$}[\atright, \solidarrow]
\mor{d}{c}{$$}[\atleft, \solidarrow]
\mor{b}{a}{$$}[\atright, \solidarrow]
\mor{b}{d}{$R$}[\atleft, \solidarrow]
\enddc\]
\noindent That is, the Hodge operators on differential forms and cochains should commute with the de Rham map. In the case of FEM, one defines $\tilde{\star}$ via the Galerkin-Hodge operator on simplicial $p$-cochains via Whitney map and $\star$ of differential forms by
\begin{equation}\label{Galerkin-Hodge}
	\langle c_1, c_2 \rangle = \int_\Omega W(c_1) \wedge \star W(c_2).
\end{equation}
This is a discrete representation of the Hodge operator, separating the metric-dependent properties in FEM from metric-free ones. Noting that $RW = I$ on cochains, we see that the Galerkin-Hodge is essentially designed to make diagrams commute. The non-degeneracy of \eqref{Galerkin-Hodge} is inherited from the $L^2$ Hodge inner product. In FEM, it results in an invertible, but not diagonal matrix, and thus, the locality of $\star$ is not strictly preserved: a manifestation of the discrete star localization problem \citep{Kotiuga2008}. The weak formulation is thus led to dictate the metric and thus the constitutive laws. The discrete problem for finding the array $\mathbf{x}$ of degrees of freedom of the problem obtains the matrix form
\begin{equation}\label{discreteGalerkin}
	\mathbf{C}^{\mathrm{T}} \mathbf{M} \mathbf{C} \mathbf{x} = \mathbf{v},
\end{equation}
where the matrix $\mathbf{M}$ is an instance of the discrete Hodge operator, the Galerkin-Hodge \cite{Bossavit1998b}. It is this Galerkin-Hodge process which allows us to roll the metric through the variational principle down to the discrete level in an essentially functorial manner. Moreover, we see here the connection to naturality related to geometric structure (see subsection~2.3.2) and how and to what extent it manifests itself via the Galerkin-Hodge process on the discrete level. Note also that the metric of the elliptic complex does not have to come from the Riemannian metric of the underlying manifold: utilization of different metrics is beneficial in tackling e.g. different constitutive laws.

Nonetheless, within the context of the category of Riemannian manifolds and isometries, Whitney forms with the Galerkin-Hodge inner product and finite element error analysis combine into a canonical procedure for discretizing inner-product-dependent properties of continuum in a geometrically compatible manner while keeping track of convergence properties with refinements of the simplicial complex. The properties of Whitney forms bring the essentially geometric nature of such a FEM procedure to the forefront \citep{KettunenTrevisan}.

\begin{myremark} Whitney forms are \emph{flat forms} arising from one-to-one correspondence with flat norm completed cochain spaces \citep{Whitney}. The exterior derivative of a flat form is flat but the Hodge star of a flat form is not: This is a fundamental consequence of discretization of metric-dependent constitutive laws of continuum physics utilizing Whitney forms. As a concise introduction to such problematics, we refer the reader to \citep{KangasAssessment}. For further considerations on flat and \emph{sharp} norm topologies and developments on the associated problematics, see Jenny Harrison's works such as \citep{Harrison}, \citep{Harrison2012}. \end{myremark} 

\subsubsection{Discrete Hodge in cochain methods}

For representing constitutive relations, cochain methods typically rely on dual complexes of so-called outer oriented chains and twisted cochains.\footnote{Some physical quantities (cochains) are naturally associated with an inner orientation of the geometric entity they are related to, while some (twisted cochains) require to be interpreted on outer oriented chains. Hence, in the primal-dual complex pair, we associate an outer oriented $n-p$-cell of the dual complex $\tilde{K}$ to each inner oriented $p$-cell of the primal complex $K$. Formally, in such a complex pair, a dual cell of a $p$-cell is such that its cofaces are the dual cells of the boundary of the $p$-cell. The outer orientation of the dual complex is induced by the inner orientation of the primal one. For example, the positive crossing direction of a 2-cell of the dual complex on 3-dimensional $\Omega$, or a 1-cell on 2-dimensional $\Omega$, is induced by the inner orientation of a 1-cell piercing it.} We will not dwell on the details of this issue here. Nonetheless, the fundamental issues of the local nature of constitutive laws and of discrete Hodge operators and how such should be defined, underlie this connection \citep{Auchmann,Hiptmair2001b,Kangas,Kangas2007,Tarhasaari,Kotiuga2008}. Another way to introduce the discrete Hodge is to express the Poincar\'e duality using a commutative cup product on cochains and define the Hodge by combining it with an inner product to avoid the dual mesh construction \citep{Wilson}.

On the matrix level, cochain methods tend to end up in a similar place as FEM: An equation of the form of \eqref{discreteGalerkin} is formed. Utilizing the same simplicial complex, the implied system matrix is the same in e.g. the cell method and FEM. However, the discrete Hodge is not the same. The essential difference in FEM and various cochain method discretizations is thus in the instantiation of the discrete Hodge -- the way metric properties and constitutive laws are incorporated into the discretization. As discussed in \citep{Tarhasaari}, equations of the form \eqref{discreteGalerkin} give rise to a circuit interpretation: The numerical approaches have similar structure to that of \emph{circuit equations} arising from circuit theory. There, the role of discrete Hodge is played by the \emph{impedance}, which connects chains to cochains. This also raises further questions on the category theoretical interpretation of discretizations and coupled discrete problems, e.g. via a monoidal categorical framework for circuit theory \citep{PassiveLinear}.\footnote{For monoidal categories, see e.g. \citep{Baez} and \citep{Coecke}.}

From this perspective, a numerical method based on an algebraic formulation of physics is not free of the questions and compatibility issues related to the differences of the continuous and the discrete, some of which we have touched upon in this paper. Error analysis is essential, and e.g. \emph{naturality} should play a similar role in the discrete descriptions as it does in the continuum descriptions. Even though there are different and differently motivated ways to approach cochain methods, such as e.g. those of Wilson and Harrison (see \citep{Wilson}, \citep{Harrison}), it all boils down to functorially preserving the continuum properties. Utilizing Whitney forms and the Galerkin-Hodge process, this is especially evident and natural.

\subsubsection{Topological properties}

Let us recap the structure we have here in terms of cohomology. The functors $L^2F^*$ and $W^*$ attach the $L^2$ de Rham complex and the Whitney complex to a manifold. Given a map $f$ between manifolds $\Omega'$ and $\Omega$, we have the induced cochain maps between the corresponding cochain complexes on $\Omega$ and $\Omega'$: $L^2F^*(f): L^2F^*(\Omega) \rightarrow L^2F^*(\Omega')$ and $W^*(f): W^*(K, \Omega) \rightarrow W^*(K', \Omega')$. These in turn induce maps on cohomology via cohomology functors. The cohomology spaces $H_s^{p}(\Omega)$ (the simplicial cohomology related to the Whitney complex) and $H_{DR}^{p}(\Omega)$ (the de Rham cohomology) are isomorphic. Moreover, between the cohomology functors of different degrees, we have the exterior derivatives (coboundaries) inducing natural transformations. So, in terms of cohomology, such a structure-preserving discretization is characterized by two commutative squares


\[\begindc{\commdiag}[50]
\obj(0,10)[a]{$H_s^{p}(\Omega)$}
\obj(20,10)[c]{$H_s^{p}(\Omega')$}
\obj(20,0)[d]{ $H_s^{p+1}(\Omega')$}
\obj(0,0)[b]{ $H_s^{p+1}(\Omega)$}
\mor{a}{c}{$H_s^{p}(f)$}[\atleft, \solidarrow]
\mor{c}{d}{$\rmd$}[\atleft, \solidarrow]
\mor{a}{b}{$\rmd$}[\atright, \solidarrow]
\mor{b}{d}{$H_s^{p+1}(f)$}[\atright, \solidarrow]

\obj(38,10)[e]{$H_{DR}^{p}(\Omega)$}
\obj(58,10)[f]{$H_{DR}^{p}(\Omega')$}
\obj(58,0)[g]{ $H_{DR}^{p+1}(\Omega')$}
\obj(38,0)[h]{ $H_{DR}^{p+1}(\Omega)$}
\mor{e}{f}{$H_{DR}^{p}(f)$}[\atleft, \solidarrow]
\mor{f}{g}{$\rmd$}[\atleft, \solidarrow]
\mor{e}{h}{$\rmd$}[\atright, \solidarrow]
\mor{h}{g}{$H_{DR}^{p+1}(f)$}[\atright, \solidarrow]
\enddc\]

\noindent for each $p \in \{0, 1,..., n-1\}$, where each of the corresponding simplicial and de Rham cohomology groups are isomorphic. Again, to be explicit, on the left $\rmd$ is the transpose of the simplicial boundary operator, and on the right it is the exterior derivative of differential forms.

\subsection{Remarks}

Our point of view emphasizes the attainment of a discrete version of the theory in question in a more general sense than just the particular solution space; the whole complex is important.\footnote{Even though the de Rham complex is the canonical example of an elliptic complex, one can consider discretizations of other complexes, too. In \citep{Arnold}, for example, the \emph{elasticity complex} and its FEM discretization are discussed.} When discretizing, Whitney forms are a very effective bridge to the cochain perspective in the context of naturality associated with differentiable structure, as well as from the perspective of essential topological properties. When it comes to naturality associated with geometric structure, we rely on the Galerkin-Hodge process, and in particular, on the Galerkin-Hodge inner product \eqref{Galerkin-Hodge} to induce an inner product on cochains. Moreover, naturality should provide us with ideas to extend these concepts to a Galerkin-based discussion of connections and curvature, relevant for e.g. generalizations of electromagnetism such as the Yang-Mills theory \citep{Garrity}. We conclude this section with a related research problem.

\begin{myrp}[Galerkin FEM and cochains: functoriality]
	The success of Whitney form Galerkin discretizations and cochain methods in CEM and continuum physics in general arises from their geometrical essence. This essence can be articulated via naturality, and often, via elliptic complexes. A central question is, how to articulate the discretization process fully functorially. Moreover, what do we ultimately mean by, or what comprises, a (structure-preserving) discretization? After all, we are searching for finite formulations of physics.
	
\end{myrp}

\section{Natural differential operators in the discrete setting}

Natural differential operators have functorial analogues on the cochains of triangulated manifolds. In order to establish this correspondence, one needs a suitable definition of a simplicial vector field. We derive formulas involving simplices and barycentric coordinates in order to define discrete vector fields and a discrete Lie derivative as a result of a discrete analogue of Cartan's magic formula. This result is articulated in Theorem~\ref{theorem1} near the end of this section.

Diffeomorphisms are tricky to talk about in the discrete setting but, with piecewise-linear (P-L) homeomorphisms, we can connect with naturality in the topological setting. Infinitesimal diffeomorphisms are vector fields and we will see that a solid definition of a discrete vector field yields a plethora of tools. Hence, we set out to develop a notion of discrete vector field which enables us to develop a discrete Lie derivative and a discrete analogue of Cartan's magic formula. To the extent that there exist discrete natural differential operations on cochains, as a major goal we would like to explore what functorial constraints ensure that they are unique.

The outline of our derivation is as follows.

\begin{enumerate}
	
	\item We give a definition of a discrete vector field in terms of a perturbation of the vertices of a simplex embedded in $\mathbb{R}^n$, using barycentric coordinates and derivatives w.r.t. Cartesian coordinates. As a check, this definition gives the proper notion of a divergence-zero vector fields as infinitesimal volume-preserving diffeomorphisms. This yields the result that discrete volume-preserving diffeomorphisms, close to the identity, are just perturbations of nodes that induce P-L transformations where the Jacobian matrices on each highest dimensional simplex have determinant 1.
	
	\item Using the definition of a discrete vector field, we show that a discrete Lie derivative follows from an analogue of Cartan's magic formula, involving the coboundary operator and this notion of infinitesimal diffeomorphism.
	
	\item We show that our discrete Lie derivative fits into a commutative diagram involving the de Rham and Whitney maps between the de Rham and Whitney complexes. As it is defined in terms of such diagrams, it is compatible with them. The Hodge-theoretic picture then follows from a variational principle and the Galerkin method.

\end{enumerate}


\subsection{Finite-element data structures in the context of naturality}

In both finite element analysis and the study of manifolds, the notion of a simplicial complex has long been used as a basic data structure enabling one to model spaces without making implicit geometrical or topological assumptions. Even though in numerical methods we tend to play with the geometrical realization of a simplicial complex, there is a rich, purely algebraic structure in data structures built from simplicial complexes, which points to more systematic use of category theory and homotopy theory in practice. 

There are three basic types of finite-element data structures: 

\begin{enumerate}
	
	\item Encoding of the abstract simplicial complex.
	
	\item Encoding of the geometric realization of an abstract simplicial complex. Contrary to what one may be inclined to think, this does not yet involve geometry per se, unless one introduces a metric on $\mathbb{R}^n$.
	
	\item The Galerkin-Hodge inner product i.e. the FE stiffness matrix. This is where the metric is introduced. All the geometry is encoded in this data structure.
	
\end{enumerate}


Given a triangulated $n$-dimensional manifold $\Omega$ (with boundary) embedded in $\mathbb{R}^n$, we can define two data structures, together characterizing the \emph{geometric realization} of the arising simplicial complex. Namely, these are the lists of globally numbered vertices
\begin{equation}\label{eqn:VertexList}
	\{1,...,m_0\}
\end{equation}
and the list of $n$-simplices given in terms of the global node numbering of the vertices
\begin{equation}\label{eqn:SimplexList}
	\{s_1, ..., s_{m_n}\},
\end{equation}
with each $s_i$ being a list of the $n+1$ vertices of the simplex. In simple terms, this is just a compact representation of the connection matrix (as per Silvester and Ferrari) \citep{SilvesterFerrari}. Indeed, \eqref{eqn:SimplexList} with each $s_i$ supplemented with $+1$ or $-1$, depending on whether the lexicographic order agrees or disagrees with the induced orientation, enocdes the data structure of an abstract simplicial complex.

The second data structure, the geometric realization of a simplicial complex, attaches to each vertex the cartesian coordinates of that vertex. That is, 
\begin{equation}\label{eqn:CartesianCoordsList}
	\{x_1, ..., x_{m_0}\},
\end{equation}
where $x_i$ is an $n$-tuple of the Cartesian coordinates attached to vertex $i$. By computing the oriented volume of each $n$-simplex $\mathrm{Vol}(s_i)$ using the Cartesian coordinates in \eqref{eqn:CartesianCoordsList} and dividing the result by its absolute value, yielding either $+1$ or $-1$,
\begin{equation}\label{eqn:AbstractSimplicialPlusMinus}
	\left\{\mathrm{Vol}(s_1)/\mid\mathrm{Vol}(s_1)\mid, ..., \mathrm{Vol}(s_{m_n})/\mid\mathrm{Vol}(s_{m_n})\mid\right\}
\end{equation}
yields the orientations of the simplices in the (abstract) simplicial complex, mentioned above. This data structure gives us a piecewise-linear structure which can be smoothed to a differentiable structure. That is, this geometric realization of the abstract simplicial complex gives us what we need to compute barycentric coordinates. By perturbing the coordinates $x_i$ we obtain the infinitesimal diffeomorphisms which we define to be the piecewise-linear vector fields.

The third data structure, the Galerkin-Hodge inner product, introduces the metric notions to FEM and thus encodes the \emph{geometry} of the problem at hand. Although there is an obstruction to having a discrete Hodge \emph{operator}, that could replicate the continuum property $(\star)^2 = (-1)^{p(n-p)} I$, (DSLP), this Galerkin-Hodge is easy to construct as
\begin{equation}\label{eqn:Galerkin-Hodge}
	\langle C_1, C_2 \rangle = \int_\Omega W(C_1) \wedge \star W(C_2),
\end{equation}
which is numerically encoded in the finite-element stiffness matrix.

In this process, we ``forget'' the geometric information present in the geometric realization and, in fact, transfer from the category of geometric realizations of abstract simplicial complexes and simplicial maps between such realizations, {\bf SGeom}, to the category of abstract simplicial complexes and simplicial maps, {\bf S}. That is, we have a \emph{forgetful functor}
\begin{equation}\label{ForgetfulFunctor}
	Abstr: {\bf SGeom} \rightarrow {\bf S}, 
\end{equation}
that takes the geometric realization to the underlying abstract simplicial complex with the operation on morphisms defined in an obvious way, i.e. simplicial maps in ${\bf SGeom}$ are taken to simplicial maps in ${\bf S}$. Functoriality of such a construction, i.e. preservation of identities and composition, is evident.

Now, having a clear distinction between an abstract simplicial complex and its geometric realization in terms of concrete finite element data structures, it is clear that many important topological computations can be done in the abstract simplicial complex without resorting to floating point arithmetic at all. We simply utilize the forgetful functor \eqref{ForgetfulFunctor} to obtain the abstract simplicial complex from its geometric realization, compute, and ``pull back'' the results to the geometric realization.

The array of vertex coordinates encoding the geometric realization of an abstract simplicial complex enables us to define piecewise-linear diffeomorphisms. The use of barycentric coordinates makes all formulas affine-invariant. Discrete vector fields are defined in terms of these P-L diffeomorphisms, and so they only require the data structure of geometric realization of a simplicial complex. Since vector fields are infinitesimal diffeomorphisms, we need only consider diffeomorphisms close to the identity. Approximating diffeomorphisms not close to the identity is irrelevant for the definition of a vector field. This would require keeping track of topological changes in the FE mesh via Pachner moves and updates to the data structure of abstract simplicial complex. 

Discrete natural differential operators require the coboundary operator and the notion of a discrete vector field to reach Cartan's magic formula. However, they do not require the graded-commutative, associative cochain algebra structure associated with wedge multiplication. Thus, discrete natural differential operators can be defined solely in terms of the first two types of data structures; the abstract simplicial complex and formulas depending on its geometric realization. 

Let us restate that our definition of a discrete vector field, and thus the discrete analogs of natural differential operators require the distinction between simplical complexes, and their geometric realizations. However, they do not involve geometry, as it is encoded in the third data structure represented by the Galerkin-Hodge inner product.

\subsection{Towards discrete Lie derivative}

Weil's concept of ``spaces of infinitely near points'' \citep{Weil1953} has led to that of \emph{Weil bundles}, which lend themselves to a characterization of all product-preserving functors on manifolds \citep{Kolar}. Inspired by the theory of jets, Weil bundles provide an algebraic framework for prolongations of manifolds within the context of naturality. This points to their utilization in the context of discretizations related to vector fields understood as infinitesimal diffeomorphisms, adducing their connection to flows and Lie derivatives and the duality between extrusion of a manifold and contraction with a vector field, as discussed by Bossavit in \citep{Bossavit03}.

\subsubsection{Formulas involving simplices and barycentric coordinates}

If $\{p_i\}^n_0$ are the $n+1$ non-coplanar points in $\mathbb{R}^n$ defining the simplex $\sigma = \langle p_o, p_1, \ldots, p_n \rangle$ and $p_k$ has coordinates $(x^k_1, x^k_2, \ldots, x^k_n)$, the $n$-volume of $\sigma$ is given by
\begin{equation}\label{eqn:simpvol}
	\mathrm{Vol}(\sigma)=\frac{1}{n!}\left| \begin{array}{ccccc}
		1 & x_{1}^{0} & x_{2}^{0} & \cdots & x_{n}^{0} \\
		1 & x_{1}^{1} & x_2^1     & \cdots & x_n^1     \\
		&  & \vdots  &  &  \\
		1 & x_{1}^{n} & x_{2}^{n} & \cdots & x_{n}^{n}  \end{array} \right|.
\end{equation}
Consider a point $p \in \sigma$ such that $\sigma$ can be subdivided into $n+1$ $n$-simplices $\sigma_k = \langle p_0,\ldots,p_{k-1},p,p_{k+1},\ldots,p_{n}\rangle$ for $0\leq k\leq n$, giving $\mathrm{Vol}(\sigma)=\sum_{k=0}^{n} \mathrm{Vol}(\sigma_{k})$. Then the $i^\mathrm{th}$ barycentric coordinate corresponding to the point $p$ is given by the fraction of the volume of $\sigma_i$ to $\sigma$ as
\begin{equation}\label{eqn:barycentric}
	\mu_i = \frac{\mathrm{Vol}(\sigma_i)}{\mathrm{Vol}(\sigma)},
\end{equation}
where
\begin{equation}\label{eqn:simpvoli}
	\mathrm{Vol}(\sigma_i)=\frac{1}{n!}\left| \begin{array}{lllll}
		1 & x_{1}^{0} & x_{2}^{0} & \cdots & x_{n}^{0} \\
		&  & \vdots &  &            \\
		1 & x_{1}^{i-1} & x_{2}^{i-1} & \cdots & x_{n}^{i-1} \\
		1 & x_{1} & x_{2} & \cdots & x_{n} \\
		1 & x_{1}^{i+1} & x_{2}^{i+1} & \cdots & x_{n}^{i+1} \\
		&  & \vdots &  &            \\
		1 & x_{1}^{n} & x_{2}^{n} & \cdots & x_{n}^{n} \end{array} \right|.
\end{equation}
Note that in the determinant of \eqref{eqn:simpvoli}, the $i^\mathrm{th}$ row is replaced by the coordinates of $p \in \sigma$. Moreover, from the definition it follows that barycentric coordinates are subject to
\begin{equation} \label{eqn:bc1}
	\sum_{k=0}^{n}\mu_{k}=1\qquad\hbox{and}\qquad
	\sum_{k=0}^{n}\mathrm{d}\mu_{k}=0.
\end{equation}

Consider a function $f\left(\left\{\mu_i\left( \left\{p_k\right\}^{n+1}_{k=1}\right)\right\}\right)$ and compute the derivative $\frac{\partial f}{\partial x^k_l}$ using the chain rule:
\begin{equation}\label{eqn:dfdx}
	\frac{\partial f}{\partial x^k_l} = \sum_{i=0}^{n} \frac{\partial f}{\partial \mu_i} \frac{\partial \mu_i}{\partial x^k_l},
\end{equation}
where
\begin{equation}\label{eqn:dmudx}
	\frac{\partial \mu_i}{\partial x^k_l} = \frac{\partial}{\partial x^k_l}\left(\frac{\mathrm{Vol}(\sigma_i)}{\mathrm{Vol}(\sigma)}\right).
\end{equation}
So, as the derivative of this rational function we obtain
\begin{equation}\label{eqn:dmudx2}
	\frac{\partial \mu_i}{\partial x^k_l} = \frac{1}{\mathrm{Vol}(\sigma)}\frac{\partial}{\partial x^k_l}\mathrm{Vol}(\sigma_i) - \frac{\mathrm{Vol}(\sigma_i)}{\left(\mathrm{Vol}(\sigma)\right)^2} \frac{\partial}{\partial x^k_l}\mathrm{Vol}(\sigma).
\end{equation}
Now, we substitute \eqref{eqn:dmudx2} to \eqref{eqn:dfdx}. This yields
\begin{equation}\label{eqn:dfdx2}
	\frac{\partial f}{\partial x^k_l} = \frac{1}{\mathrm{Vol}(\sigma)} \sum_{i=0}^{n} \frac{\partial f}{\partial \mu_i} 	\left[\frac{\partial}{\partial x^k_l} \mathrm{Vol}(\sigma_i) - \frac{\mathrm{Vol}(\sigma_i)}{\mathrm{Vol}(\sigma)} \frac{\partial}{\partial x^k_l} \mathrm{Vol}(\sigma) \right]
\end{equation}	
or more succinctly
\begin{equation}\label{eqn:dfdx3}
	\frac{\partial f}{\partial x^k_l} = \frac{1}{\mathrm{Vol}(\sigma)} \sum_{i=0}^{n} \frac{\partial f}{\partial \mu_i} 	\left[\frac{\partial}{\partial x^k_l} \mathrm{Vol}(\sigma_i) - \mu_i \frac{\partial}{\partial x^k_l} \mathrm{Vol}(\sigma) \right].
\end{equation}
The partial derivatives of the volumes of the simplices with respect to $x^k_l$ in \eqref{eqn:dfdx3} are given as
\begin{equation}\label{eqn:dvolsigmaidx}
	\frac{\partial}{\partial x^k_l} \mathrm{Vol}(\sigma_i) = \frac{1}{n!}\left| \begin{array}{llllllllll}
		1 & x_{1}^{0} & x_{2}^{0} & \cdots & & & & & x_{n}^{0} \\
		1 & x_{1}^{1} & x_{2}^{1} & \cdots & & & & & x_{n}^{1}  \\
		&  & \vdots &  & & & &  &         \\
		1 & x_{1}^{k-1} & x_{2}^{k-1} & \cdots & & & & & x_{n}^{k-1} \\
		0 & 0 & 0 & \cdots & 0 & 1 & 0 & \cdots & 0 \\
		1 & x_{1}^{k+1} & x_{2}^{k+1} & \cdots & & &  & & x_{n}^{k+1} \\
		&  & \vdots &  &  & & & &          \\
		1 & x_{1}^{n} & x_{2}^{n} & \cdots & & & & & x_{n}^{n} \end{array} \right|,
\end{equation}
for $i \neq k$. Note that row $k$ consists of zeros except for the $(l+1)^\mathrm{st}$ entry which is $1$. Moreover, note that the $i^\mathrm{th}$ row consists of $(1, x_1, x_2, \cdots, x_n)$. Thus, simplifying the determinant, we can write
\begin{equation}\label{eqn:dvolsigmaidx2}
	\frac{\partial}{\partial x^k_l} \mathrm{Vol}(\sigma_i) = \frac{(-1)^{(k+l)}}{n!}\left| \begin{array}{llllllllll}
		1 & x_{1}^{0} & x_{2}^{0} & \cdots & x_{l-1}^{0} & x_{l+1}^{0} & \cdots & & x_{n}^{0} \\
		1 & x_{1}^{1} & x_{2}^{1} & \cdots & x_{l-1}^{1} & x_{l+1}^{1} & \cdots & & x_{n}^{1} \\ 
		&  & \vdots &  & & & &  &         \\
		1 & x_{1}^{k-1} & x_{2}^{k-1} & \cdots &   x_{l-1}^{k-1} & x_{l+1}^{k-1} & \cdots & & x_{n}^{k-1} \\
		1 & x_{1}^{k+1} & x_{2}^{k+1} & \cdots & x_{l-1}^{k+1} & x_{l+1}^{k+1} & \cdots & & x_{n}^{k+1} \\
		&  & \vdots &  &  & & & &          \\
		1 & x_{1}^{n} & x_{2}^{n} & \cdots & x_{l-1}^{n} & x_{l+1}^{n} & \cdots & & x_{n}^{n} \end{array} \right|.
\end{equation}
We can interpret this as the volume of an $(n-1)$-simplex formed by projecting the face of $\sigma$ opposite to $p_k$ onto the $(n-1)$-dimensional hyperplane obtained by setting $x_l = 0$.

The partial derivative $\frac{\partial}{\partial x^k_l} \mathrm{Vol}(\sigma)$ is computed similarly as \eqref{eqn:dvolsigmaidx2} but without the $i \neq k$ restriction. Substituting these partial derivatives to \eqref{eqn:dfdx3} gives us an explicit formula for $\frac{\partial f}{\partial x^k_l}$.

\subsubsection{Discrete vector fields via 1-parameter deformation of vertices}

Next we define discrete vector fields in a manner compatible with the de Rham complex. 

We start with a warm-up problem that hints at the definition of a discrete vector field. Expressing the warm-up problem in terms of barycentric coordinates inspires the next step. We need to consider a parameterized deformation of a simplicial complex, where all the vertices are perturbed as a function of $t$, while maintaining the affine geometry of the barycentric coordinates. In order to obtain a 1-parameter-family of complexes deforming the original complex, let $\intprod$ denote contraction, and on a particular simplex, consider $f\left(\left\{\mu_i\left(\mathrm{vertices}(t)\right) \right\}\right)$. So, by the chain rule we can write
\begin{equation}\label{eqn:dfdtpervertex}
	\frac{\rmd f}{\rmd t} = \sum_{i,k,l=0}^{n}\frac{\partial f}{\partial \mu_i}\frac{\partial \mu_i}{\partial x^k_l}\frac{\rmd x^k_l}{\rmd t}
\end{equation}
which, using \eqref{eqn:dmudx2}, becomes
\begin{equation}\label{eqn:dfdtpervertex2}
	\frac{\rmd f}{\rmd t} = \sum_{i,k,l=0}^{n} \frac{\partial f}{\partial \mu_i}\left[ \frac{1}{\mathrm{Vol}(\sigma)} \left(\frac{\partial \mathrm{Vol}(\sigma_i)}{\partial x^k_l} - \mu_i \frac{\partial\mathrm{Vol}(\sigma)}{\partial x^k_l} \right) \right]  \frac{\rmd x^k_l}{\rmd t}.
\end{equation}
We want to interpret this as $\frac{\rmd f}{\rmd t} = v \intprod \rmd f$. There are two sensible ways to go about this. First, with $\{\mu_i\}^k_{i=0}$ as coordinates, we have:
\begin{mylemma}[Discrete contraction 1]
\begin{equation}
	\frac{\rmd f}{\rmd t} = \sum_{i=0}^{n} \frac{\partial f}{\partial  \mu_i} \frac{\rmd \mu_i}{\rmd t} =  v \intprod \rmd f,
\end{equation}
where 
\begin{equation}\label{eqn:df(mu_i)}
	\rmd f = \sum_{i=0}^{n} \frac{\partial f}{\partial \mu_i} \rmd \mu_i, 
\end{equation}	
\begin{equation}
	v = \sum_{i=0}^{n} \frac{\rmd \mu_i}{\rmd t} \frac{\partial}{\partial \mu_i},
\end{equation}		
and $\frac{\partial}{\partial \mu_i} \intprod \rmd \mu_j$ is subject to
\begin{equation}
	\sum_{i=0}^{n} \mu_i = 1, \quad \mathrm{or} \quad 	\sum_{i=0}^{n} \rmd \mu_i = 0.
\end{equation}
	
\end{mylemma}
\noindent This is quite intrinsic. However, this obscures the connection to a P-L deformation of the mesh. Alternatively, if we explicitly use the embedding of the mesh into $\mathbb{R}^n$, and use formulas \eqref{eqn:simpvoli}--\eqref{eqn:df(mu_i)}, we consider the deformation of the vertices which induce $\{\frac{\partial {\mu_i}}{\partial t}\}^k_{i=0}$ so that
\begin{mylemma}[Discrete contraction 2]
\begin{equation}
	\frac{\rmd f}{\rmd t} = \sum_{k,l=0}^{n} \frac{\partial f}{\partial  x^k_l} \frac{\rmd x^k_l}{\rmd t} =  v \intprod \rmd f, 
\end{equation}
where
\begin{equation}
	\rmd f = \sum_{k,l=0}^{n} \frac{\partial f}{\partial x^k_l} \rmd x^k_l, 
\end{equation}	
\begin{equation}
	v = \sum_{k,l=0}^{n} \frac{\rmd x^k_l}{\rmd t} \frac{\partial}{\partial x^k_l},
\end{equation}	
subject to
\begin{equation}
	\frac{\partial}{\partial x^i_j} \intprod \rmd x^k_l = \delta_{ik}\delta_{jl}.
\end{equation}
\end{mylemma}

\subsection{A discrete version of Cartan's magic formula}

We have introduced the notion of a discrete vector field as an infinitesimal diffeomorphism, and the idea that contractions of cochains with discrete vector fields result from the chain rule. But the obvious candidate for Cartan's magic formula still remains to be found. 

Let us re-examine how the Lie derivative arises on differential forms from first principles, and mimic its derivation in the cochain context. This leads us to the notion of a discrete version of Cartan's magic formula. For relevant background, see \citep{Frankel}, chapter 4.

But first, we turn our attention to the Whitney and de Rham maps, so we can have the whole construction follow from existing definitions.

\subsubsection{Dupont's chain homotopy operator $K_\mathrm{c}$}

Let us recall the de Rham ($R$) and Whitney ($W$) maps defined in \eqref{deRhamMap} and \eqref{WhitneyMap}, respectively. 

For $K$, a triangulation of a manifold $\Omega$, we have 
\begin{equation}
	R:F^p(\Omega)\rightarrow C^p(K), 
\end{equation}
\begin{equation}
	W:C^p(K)\rightarrow F^p(\Omega). 
\end{equation}
Then, the following holds
\begin{equation}
	RW=I:C^p(K)\rightarrow C^p(K)
\end{equation}
\begin{equation}
	WR:F^p(\Omega)\rightarrow F^p(\Omega).  
\end{equation}
So, $RW$ is indeed the indentity map, but for $WR$, which is a projection onto the image of $W$ in the de Rham complex, we have the \emph{deformation retract}\footnote{A deformation retract is a homotopy equivalence with one of its homotopies being the identity.} property of \emph{Dupont's chain homotopy operator} $K_\mathrm{c}$ \citep{Dupont}
\begin{equation}
	I-WR=\rmd K_\mathrm{c}+K_\mathrm{c} \rmd:F^p(\Omega)\rightarrow F^p(\Omega). 
\end{equation}
Or, in terms of diagrams:

\[\begindc{\commdiag}[50]
\obj(0,40)[null1]{$C^p(K)$}
\obj(0,20)[a]{$C^p(K)$}
\obj(30,20)[b]{$F^p(\Omega)$}
\obj(30,40)[null2]{$F^p(\Omega)$}
\mor{a}{b}{$W$}
\mor{b}{null1}{$R$}[\atright, \solidarrow]
\mor{b}{null2}{$I - \rmd K_\mathrm{c} - K_\mathrm{c} \rmd$}[\atright, \solidarrow]
\mor{a}{null1}{$I$}
\mor{null1}{null2}{$W$}

\enddc\].

Now, let us explicitly define the operator $K_\mathrm{c}$ that appears above. A chain homotopy is an algebraic construction ensuring that a homotopy induces an isomoprhism between homology groups. As a familiar example, consider the mapping cone construction, which shows that a contractible space has the same homology groups as a point (Poincare Lemma). In our context, Dupont's chain homotopy operator ensures that there is a chain homotopy between $WR$ and the identity mapping $I$. 

Dupont's chain homotopy operator $K_\mathrm{c}$ is a mapping
\begin{equation}
	K_\mathrm{c}:F^p(\Omega)\rightarrow F^{(p-1)}(\Omega) 
\end{equation}
with properties
\begin{equation}
	K_\mathrm{c}W=0: C^p(K)\rightarrow F^{(p-1)}(\Omega), 
\end{equation}
\begin{equation}
	RK_\mathrm{c}=0:F^p(\Omega)\rightarrow C^{(p-1)}(K), 
\end{equation}
\begin{equation}
	\rmd K_\mathrm{c}+K_\mathrm{c} \rmd =I-WR:F^p(\Omega)\rightarrow F^p(\Omega), 
\end{equation}
\begin{equation}
	K_\mathrm{c}^2=0:F^p(\Omega)\rightarrow F^{(p-2)}(\Omega).
\end{equation}
We may write $F^*(\Omega)$ as
\begin{equation}\label{eqn:dirsum1}
	F^*(\Omega)= \mathrm{Im}(W)\bigoplus \mathrm{Ker}(R).
\end{equation}
Now, $K_\mathrm{c}: \mathrm{Ker}(R)\rightarrow \mathrm{Ker}(R)$,
$WR = I: \mathrm{Im}(W) \rightarrow \mathrm{Im}(W)$, and since $RK_\mathrm{c}=0$, and $K_\mathrm{c}=\rmd^{-1}$ on $\mathrm{Ker}(R)$, we have
\begin{equation}\label{eqn:dirsum2}
	\mathrm{Ker}(R)= (\mathrm{Im}(\rmd)\cap \mathrm{Ker}(R))\bigoplus \mathrm{Im}(K_\mathrm{c}).
\end{equation}
Combining \eqref{eqn:dirsum1} and \eqref{eqn:dirsum2}, we obtain a Hodge-like decomposition:
\begin{equation}
	F^*(\Omega)= \mathrm{Im}(W)\bigoplus (\mathrm{Im}(\rmd)\cap \mathrm{Ker}(R))\bigoplus \mathrm{Im}(K_\mathrm{c}).
\end{equation}

\begin{mydef}
Define $\rmd_\mathrm{Whitney} = W\partial^\mathrm{T}R$,  where $\partial^\mathrm{T}$ is the coboundary operator. \\
\end{mydef}
\noindent Then on $\mathrm{Im}(W)\bigoplus (\mathrm{Im}(\rmd)\cap \mathrm{Ker}(R))\bigoplus \mathrm{Im}(K_\mathrm{c})$, we may write the exterior derivative of the de Rham complex as a "block matrix"
$$
\rmd=\left[\begin{array}{ccc}
	\rmd_\mathrm{Whitney} & 0 & 0 \\
	0 & 0 & \rmd_\mathrm{uv} \\
	0 & 0& 0
\end{array}\right]
$$
and Dupont's chain homotopy operator as\footnote{The symbols $\mathrm{uv}$ here refer to ultraviolet in the sense of lattice gauge theory.}
$$
K_\mathrm{c}=\left[\begin{array}{ccc}
	0 & 0 & 0 \\
	0 & 0 & 0 \\
	0 & \rmd_\mathrm{uv}^{-1} & 0
\end{array}\right].
$$

\noindent This allows us to summarize the action of $\rmd_\mathrm{Whitney}$ on our Hodge-like decomposition by

\[\begindc{\commdiag}[50]
\obj(0,20)[a]{$C^{p}(K)$}
\obj(0,10)[b]{$F^{p}(\Omega)$}
\obj(25,10)[c]{$F^{p+1}(\Omega)$}
\obj(25,0)[d]{$C^{p+1}(K)$}
\obj(0,0)[e]{$C^{p}(K)$}
\obj(25,20)[f]{$F^{p}(\Omega)$}

\mor{a}{f}{$W$}[\atleft, \solidarrow]
\mor{b}{a}{$R$}[\atleft, \solidarrow]
\mor{b}{e}{$R$}[\atright, \solidarrow]
\mor{f}{c}{$\rmd$}[\atleft, \solidarrow]
\mor{d}{c}{$W$}[\atleft, \solidarrow]
\mor{e}{d}{$\partial^\mathrm{T}$}[\atright, \solidarrow]
\mor{b}{c}{$\rmd_{\mathrm{Whitney}}$}[\atleft, \solidarrow]
\enddc\].	

\subsubsection{Discrete version of Cartan's magic formula}

Let us first consider the realm of differential forms. Consider a domain $D$ under the action of a 1-parameter family of diffeomorphisms $\varphi(x,t)$ so that $D$ becomes $D(t)$. Then
\begin{equation}\label{diffInt}
	I = \int_D \omega \rightsquigarrow I(t) = \int_{D(t)} \omega(t)	
\end{equation}	
and
\begin{equation}\label{diffIntTimeDerivative}
	\frac{\rmd I}{\rmd t} = \lim_{\Delta t \to 0}  \frac{1}{\Delta t}\left[ \frac{I(t+\Delta t) - I(t)}{\Delta t} \right].
\end{equation}	
Or, adding zero,
\begin{equation}\label{diffIntTimeDerivative2}
	\frac{\rmd I}{\rmd t} = \lim_{\Delta t \to 0} \frac{1}{\Delta t}\left[  \int_{D(t + \Delta t)} \omega(t + \Delta t) -  \int_{D(t + \Delta t)} \omega(t) +  \int_{D(t + \Delta t)} \omega(t) -  \int_{D(t)} \omega(t)	\right]
\end{equation}
or
\begin{equation}\label{diffIntTimeDerivative3}
	\frac{\rmd I}{\rmd t} = \int_{\partial D} v \intprod \omega + \int_D \frac{\rmd \omega}{\rmd t} = \int_{D} \rmd (v \intprod \omega ) + v \intprod \rmd \omega ,
\end{equation}
where $v|_{x} = \frac{\partial \varphi} {\partial t}|_{x,0}$.


Let us mimic the above calculation leading up to \eqref{diffIntTimeDerivative3} in the realm of cochains. Specifically, let 
\begin{itemize}
	
	\item $\sigma_n$ be a highest-dimensional simplex in a triangulation.
	
	\item $\sigma_n(t)$ be a 1-parameter family of such simplices under the action of a piecewise-linear affine transformation. 
	
	\item $C_n(t)$ be an $n$-cochain on $\sigma_n(t)$ which is $R(W(C_n(t))) = R(\omega(t))$, where $\omega(t) = W(C_n(t))$ and, as usual, $W$ is the Whitney map, and $R$ is the de Rham map, that satisfy $RW = I$.
\end{itemize}

Now consider $I(t) = \int_{\sigma_n(t)} C_n = \int_{\sigma_n(t)} R(\omega(t))$. However, we are not actually computing any integrals, so let us write this as
\begin{equation}\label{I(t)}
	I(t) = \left\langle R(\omega(t)), \sigma_n(t) \right\rangle
\end{equation}
so that
\begin{equation}\label{dIdtcochain}
	\frac{\rmd I}{\rmd t} = \left\langle R\left(\frac{\partial\omega}{\partial t}\right), \sigma_n \right\rangle + \left\langle R(\omega), \frac{\partial \sigma_n}{\partial t} \right\rangle.
\end{equation}

So what exactly is going on? Let us think in terms of the familiar de Rham complex again. 
Start with Cartan's magic formula for the Lie derivative $\mathcal{L}_X$.
\begin{equation}\label{deRhamLx}
	\mathcal{L}_X{\omega} = X \intprod \rmd \omega + \rmd ( X \intprod \omega): F^p \rightarrow F^p.
\end{equation}	

Let $\omega = W(c)$ so that 
\begin{equation}\label{deRhamLx2}
	\mathcal{L}_X(W(c)) = X \intprod \rmd W(c) + \rmd(X \intprod W(c)).
\end{equation}	
However, using the commutative diagram $\rmd W = W \partial^\mathrm{T}$ we have
\begin{equation}\label{deRhamLx3}
	\mathcal{L}_X(W(c)) = X \intprod W(\partial^\mathrm{T} c) + \rmd (X \intprod W(c))
\end{equation}	
and in the realm of cochains we then have
\begin{equation}\label{deRhamLx4}
	R(\mathcal{L}_X(W(c))) = R( X \intprod W(\partial^\mathrm{T} c) ) + R( \rmd (X \intprod W(c))  ).
\end{equation}
Let $\omega \in \mathrm{Im}(W)$, i.e., $\omega = W(c)$ for some $c \in C^p(K)$. Given Cartan's magic formula \eqref{deRhamLx2},
we let $c=R(\omega)$. This yields
\begin{equation}\label{eqn:LxWR}
	\mathcal{L}_X WR(\omega) = X  \intprod (\rmd W R (\omega) ) + \rmd (X \intprod W R(\omega)).
\end{equation}
In the context of Whitney-form-based discretizations, i.e. the maps $W$ and $R$, we can thus write the following commutative diagrams for the Lie derivative, based on \eqref{deRhamLx2}, \eqref{deRhamLx4}, and \eqref{eqn:LxWR}:

\[\begindc{\commdiag}[50]
\obj(0,10)[a]{$C^{p}(K)$}
\obj(25,10)[c]{$F^{p}(\Omega)$}
\obj(25,0)[d]{$F^{p}(\Omega)$}
\obj(0,0)[b]{$F^{p}(\Omega)$}
\mor{a}{c}{$W$}[\atleft, \solidarrow]
\mor{c}{d}{$\mathcal{L}_X$}[\atleft, \solidarrow]
\mor{b}{a}{$R$}[\atleft, \solidarrow]
\mor{b}{d}{$X\intprod\rmd W R + \rmd X \intprod W R$}[\atright, \solidarrow]
\enddc\]

\[\begindc{\commdiag}[50]
\obj(0,10)[a]{$C^{p}(K)$}
\obj(25,10)[c]{$F^{p}(\Omega)$}
\obj(25,0)[d]{$F^{p}(\Omega)$}
\obj(0,0)[b]{$C^{p}(K)$}
\mor{a}{c}{$W$}[\atleft, \solidarrow]
\mor{c}{d}{$\mathcal{L}_X$}[\atleft, \solidarrow]
\mor{b}{a}{$RX\intprod W \partial^\mathrm{T} + R \rmd X \intprod W$}[\atleft, \solidarrow]
\mor{d}{b}{$R$}[\atright, \solidarrow]
\enddc\].


Having the action of the vector field on the vertices of the simplicial mesh defined in terms of the barycentric coordinates and their derivatives, and the definition of $\rmd_\mathrm{Whitney} = W\partial^\mathrm{T}R$, the above suggests a canonical way to define the discrete Lie derivative. The following theorem is the main original result of our paper.
\begin{mytheorem}[Discrete Cartan's magic formula]\label{theorem1}
\begin{equation*}\label{eqn:LxWR2}
	\mathcal{L}_x WR(\omega) = X  \intprod ( \rmd_\mathrm{Whitney} \omega ) + \rmd (X \intprod W R(\omega)).
\end{equation*}
\end{mytheorem}
\noindent {\bf Proof.} From the commutative diagram, $\rmd W = W \partial^\mathrm{T}$, we have $\rmd WR ( \omega ) = W \partial^\mathrm{T} R(\omega) = \rmd_\mathrm{Whitney} \omega$. $\square$ \\

Now, let $c(t)$ be a parametrized simplex.
\begin{equation}\label{eqn:LxBackToCalc} 
	\int_{c(t)} \mathcal{L}_X(\omega) = \frac{\rmd}{\rmd t} \int_{c(t)} \omega(t) = \int_{\partial c(t)} X \intprod \omega + \int_{c (t)} X \intprod \rmd \omega = \int_{c(t)} \rmd (X \intprod \omega) + X \intprod (\rmd \omega).
\end{equation}
Since $\omega \in \mathrm{Im}(W)$, $\omega = WR(\omega)$. Hence, combined with our definition of a discrete vector field, we have a discrete cersion of Cartan's magic formula that is compatible with $R$ and $W$. Concretely, we can write
\begin{equation}\label{eqn:discreteCartan} 
	\int_{c(t)} \mathcal{L}_X(\omega) = \int_{\partial c(t)} X \intprod WR(\omega) + \int_{c (t)} X \intprod \rmd WR(\omega).
\end{equation}
This is how we make sense of equation \eqref{dIdtcochain}.



\subsection{Remark on combinatorial cochain algebras}

Homotopy theory and category theory are tied to FEM and cochain methods also via the work of D. Sullivan and S. O. Wilson. The combinatorial analogue of the exterior product of differential forms, the cup product of cochains defined in a simplicial complex via $\tau_1 \cup \tau_2 = R(W\tau_1 \wedge W\tau_2)$, is a product in a $C_\infty$-algebra, which is (an almost) commutative ring with commutativity relaxed up to higher homotopy. Wilson shows in \citep{Wilson}, that in the continuum limit, all the higher homotopies converge to zero and the algebra converges to the commutative and associative algebra of differential forms given by the exterior product. However, again, in the discrete side, we are bound to work \emph{up to homotopy}. For the local construction of $\infty$-structures used, see e.g. the appendix to \citep{Sullivan} by D. Sullivan.

We end this section with two open research problems.

\begin{myrp}[Categories and data structures for discrete solvers]
	
	The rich algebraic structure of abstract simplicial complexes, separate from its geometrical realization, points to richer use of category theory in practice of programming modelling software. Utilizing functors from the category of abstract simplicial complexes to categories of algebraic topology, computations can be done without recourse to floating point arithmetic.
	
\end{myrp}

\begin{myrp}[Higher categories, type theory and discrete solvers]
	
	There is an intimate relationship between homotopy theory and higher-category theory \citep{Riehl}. In particular, there is a connection to homology computation \citep{Kotiuga1989b}, \citep{Pellikka} and combinatorial cochain algebras \citep{Wilson}. This is further tied to computational interpretations through homotopy type theory, type theory in general \citep{HoTT}, \citep{Angiuli}, \citep{Angiuli2017}, \citep{Shulman} and functional programming \citep{Milewski}. This is a vast, unexplored world in terms of computational physics software. 
	
\end{myrp}

\section{Conclusions}

Geometric methods have arguably been the most successful ones in computational physics of continuum, yielding structure-preserving discretizations. A highly practical setting for geometric methods is that of a natural elliptic complex. Elliptic complexes arise naturally in computational physics of continuum in compact domains, give rise to a rich, generalized Hodge theory and allow us to formulate computations in terms of the symbol sequence, consisting of finite-dimensional spaces and matrices. With well-posed boundary conditions, this becomes a Fredholm complex. Moreover, regarding programming modelling software, such a cochain complex gives rise to natural data structures for discrete solvers. From the point of view of an engineer programming such software, a general Hilbert complex is often fruitless. Furthermore, in the discrete setting, the functorial relationship between abstract simplicial complexes and their geometric realizations on a Riemannian manifold enable computations without recourse to floating-point arithmetic for some phases of the computational process. 

In this paper, we have discussed this geometric essence in the context of the finite element method and cochains. Evidently, the essentially category theoretical concept of naturality comes up at different levels of abstraction when discussing computational physics of continuum, and aspects of naturality shall be preserved in a discretization. This is the essence of geometric methods: Geometric means natural in this context. A very comprehensive treatment of this viewpoint beyond the scope of this paper can be found in \citep{Kolar}. There exist different types of naturality related to different mathematical structures, and they are important in different settings. However, a unification through category theory is possible. 

Coming from continuum, one wants to map the metric aspects of the problem in a functorial manner to the discrete setting. The Galerkin-Hodge process with Whitney forms presents itself as a canonical (i.e. functorial in the given categorical setting) way to achieve this goal. Throughout, we have emphasized the functorial nature behind structure-preserving discretizations, and especially, albeit without explicitly defining such a functor, the functorial nature of the Galerkin-Hodge process as well as the role of naturality. At the same time, this functoriality brings the formally similar structures and final goals of the Whitney form FEM approach and cochain method approach to the forefront. In the context of discretizations, Whitney's interpolation formula is uniquely characterized as the structure-preserving discretization of natural differential operators on a manifold, compatible with the structure of the complex arising from the triangulation.

As the main mathematical result of this paper (Theorem 1), we derived a discrete version of Cartan's magic formula from a definition of a discrete vector field induced by a 1-parameter deformation of the vertices of a finite-element mesh. Discrete contractions (Lemma 1 and Lemma 2) and defining $\rmd_\mathrm{Whitney}$ (Definition 1) were instrumental for this. More than that, our approach yields a conceptual tool to discuss other notions of naturality important for computational physics: Hermitian geometry, complex geometry, and connections on bundles, to name a few, should be investigated in this context of variational principles.

Looking at computational physics from an electrical engineering perspective suggests viewing several concepts of modern mathematics, not so well-established in computational physics literature, as ``missed opportunities'' in the field. These connect through category theory to form a set of research problems we set out to pursue. These are the research problems RP 1-6 distributed in the present text. Observing the categorical structures and structural similarities in discrete solvers yields us tools for abstract and effective communication between scientists, engineers and computers. We see this as an opportunity. 

\section*{Acknowledgments}

A significant part of the research leading to this paper was originally conducted during Valtteri Lahtinen's (V.L.) research visit at the ECE department of Boston University (BU) during 2017 and 2018. V.L. would like to thank his colleagues at BU for their hospitality during his stay. We would like to also acknowledge that Dr. Antti Stenvall had an instrumental influence for the very early stages of this paper. Moreover, the authors would like to especially thank Alain Bossavit, Andr\'e Nicolet and Jim Stasheff for their comments and criticism concerning the manuscript.

\section*{Funding and/or Conflicts of Interests/Competing Interests}

The authors declare that they have no conflict of interest.

\pagebreak

\end{document}